\documentclass[11pt, reqno]{amsart}
\usepackage{amsmath, amssymb}
\usepackage[all]{xypic}
\usepackage{psfig}
\newcommand{\eeq}{\end{equation}}

\newtheorem{theorem}{Theorem}[section]

\newtheorem{defi}[theorem]{Definition}
\newtheorem{lemma}[theorem]{Lemma}

\def\slfrac#1#2{\hbox{\kern.1em %
 \raise.5ex\hbox{\the\scriptfont0 #1}\kern-.11em %
 /\kern-.15em\lower.25ex\hbox{\the\scriptfont0 #2}}}

\newcommand{\fot}{\frac{1}{2}}
\newcommand{\eqn}[1]{(\ref{#1})}
\newcommand{\hsp}{\hspace*{\parindent}}

\newcommand{\beql}[1]{\begin{equation}\label{#1}}
\newcommand{\bsq}{{\vrule height .9ex width .8ex depth -.1ex }}

\font\phvr=phvr at 11pt
\newcommand\he[1]{\mbox{\phvr  #1}}

\newcommand{\tO}{[\tau_1]}
\newcommand{\tT}{[\tau_2]}

\newcommand{\pt}{\partial}

\newcommand{\ep}{\epsilon}

\newcommand{\Om}{\Omega}

\newcommand{\ZZ}{{\mathbb Z}}
\newcommand{\RR}{{\mathbb R}}

\newcommand{\QQ}{{\mathbb Q}}
\newcommand{\CC}{{\mathbb C}}

\newcommand{\bR}{{\bf R}}

\newcommand{\bx}{{\bf x}}

\newcommand{\hG}{{\he{G}}}

\newcommand{\ot}{{\phi_{[\tau]}}}

\newcommand{\sA}{{\mathcal A}}
\newcommand{\sB}{{\mathcal B}}

\newcommand{\sG}{{\mathcal G}}

\newcommand{\sM}{{\mathcal M}}

\newcommand{\sO}{{\mathcal U}}

\newcommand{\sQ}{{\mathcal Q}}
\newcommand{\sR}{{\mathcal R}}

\newcommand{\sU}{{\mathcal U}}
\newcommand{\sV}{{\mathcal V}}
\newcommand{\sW}{{\mathcal W}}
\newcommand{\sX}{{\mathcal X}}
\newcommand{\sY}{{\mathcal Y}}

\newcommand{\rC}{{\rm C}}

\newcommand{\tsM}{{\widetilde{\mathcal M}}}
\newcommand{\tOm}{\tilde{\Omega}}

\newcommand{\bth}{{\mathbf \theta}}

\title{The Lerch Zeta Function
II. Analytic Continuation}
\subjclass[2000]{Primary: 11M35}
\keywords{functional equation, Hurwitz zeta function, Lerch zeta function}

\author{Jeffrey C. Lagarias}
\thanks{The work of the first author was supported by NSF grants DMS-0500555
and DMS-0801029
and  the second author by NSF grants DMS-0457574 and
DMS-0801096.}
\address{Department of Mathematics, University of Michigan,
Ann Arbor, MI 48109-1043,USA}
\email{lagarias@umich.edu}
\author{Wen-Ching Winnie Li}
\address{Department of Mathematics, Pennsylvania State University,
University Park, PA 16802-8401, USA}
\email{ wli@math.psu.edu}

\date{June 5, 2010}
\begin{document}
\begin{abstract}
This is the second of a series of four papers that study algebraic
and analytic structures associated with the Lerch zeta function.
The  Lerch zeta function
$
\zeta (s,a,c) := \sum_{n=0}^\infty \frac{e^{2 \pi in a}}{ (n+c)^{s}}
$
was introduced by Lipschitz  in 1857,
and is  named after Lerch, who showed in 1887  that it satisfied
a functional equation. Here we analytically continue
$\zeta (s, a, c)$ as a function of three complex variables.
We show that it is well-defined as a multivalued  function
on the manifold
$\sM:= \{ (s, a, c) \in \CC \times ( \CC \smallsetminus \ZZ ) \times
( \CC \smallsetminus \ZZ ) \},$ and that
 this analytic continuation becomes single-valued on the maximal abelian cover of $\sM$.
We compute the monodromy functions describing the multivalued
nature of this function on $\sM$, and determine various of
its properties.
\end{abstract}

\maketitle
\tableofcontents
\setlength{\baselineskip}{1.2\baselineskip}

%

\section{Introduction}
\hsp
The {\em Lerch zeta function} is defined by the Dirichlet series
\beql{101}
\zeta (s,a,c) := \sum_{n=0}^\infty \frac{e^{2 \pi in a}}{ (n+c)^{s}}.
\eeq
It already appeared in a paper of
Lipschitz \cite{Li57} in 1857, but
is named after Lerch \cite{Le87}, who showed in 1887
that it satisfies a functional equation. This functional
equation, called {\em Lerch's transformation formula}, states that
\beql{102}
\zeta (1-s, a,c) =
(2 \pi )^{-s} \Gamma (s) \left\{
e^{\frac{\pi is}{2}} e^{-2 \pi iac} \zeta (s, -c,a)~
+~ e^{- \frac{\pi is}{2}} e^{2 \pi ic (1-a)} \zeta (s,c,1-a) \right\} ~,
\eeq
and is valid
for complex $a$ with $\Im(a) \ge 0$ and real $c$ with
$0 < c < 1$, cf. Erdelyi \cite[p. 29]{Er53}.
(For $\Im (a) > 0 $ and real $ c > 0 $, the Dirichlet
series \eqn{101} is an entire function of $s$, because its coefficients
are rapidly decreasing.)
In part I we derived two
 symmetrized four term functional equations, as
given by Weil \cite{We76}, which are valid for real $a$ and $c$, and for all
$s \in \CC.$  These apply to the two  functions 
$$
L^{\pm}(s, a, c) := \zeta(s, a, c) \pm e^{-2\pi i  a} \zeta(s, 1-a, 1- c),
$$
see Theorem~\ref{th11} below.
Special cases of the Lerch zeta function
include $a=0$ which gives the
{\em Hurwitz zeta function}
\beql{103}
\zeta (s,c) = \sum_{n=0}^\infty \frac{1}{(n+c)^s} ~,
\eeq
named after Hurwitz \cite{Hu82},
and the case $c=1$, which gives 
\beql{104}
e^{-2\pi a}  F(s,a) =e^{-2 \pi ia} \left(  \sum_{n=1}^\infty \frac{e^{2 \pi i na}}{n^s} \right),
\eeq
where $F(s, a)$
  is the {\em periodic zeta function}, cf. Apostol \cite[p. 257]{Ap76}. 
The Riemann zeta function
occurs when $a = 0$ and $c = 1$.  These four term functional equations,
stated in Theorem~\ref{th11}below, 
can be used to derive the known functional equations of the Hurwitz 
and periodic zeta functions.

This paper is the second in a series of four papers that study
alebraic and analytic structures associated with the Lerch zeta function.
Our object here is to obtain the analytic continuation of the Lerch
zeta function as a function of three complex variables.
Our main  result  analytically continues the function
 to an essentially  maximal domain of holomorphy, viewing it as a multivalued function
 of three variables. Then we determine properties of the extension to this
 domain, noting particularly the differential-difference equations
 and differential equations which
 these functions satisfy, described below and in more detail in \S2.

First, we  analytically continue the  Lerch zeta function
to a single-valued function on the universal cover $\tsM$ of
 the non-simply connected
manifold
\beql{106}
 \sM := \{ (s, a, c)  \in \CC \times (\CC
\smallsetminus \ZZ) \times (\CC \smallsetminus \ZZ)\}. \eeq The
resulting function may alternatively be viewed as a multivalued
function on $\sM$, and we compute the monodromy functions describing
the multivaluedness. We show that the  resulting function becomes
single valued on the maximal abelian cover $\widetilde{\sM}^{ab}$ of $\sM$; see \S2.
The manifold $\sM$ proves convenient to study because
this  is the largest manifold on  which the two  functions $L^{\pm}(s,a,c)$
appearing in the functional
equations  are well-defined as multivalued functions.
That is, the  manifold $\sM$  has the symmetry that the  map 
$\phi_{\bR}: \sM \to \sM$ defined by
$$\phi_{\bR} (s, a,c)= (1-s, 1-c, a)$$
 defines an automorphism
of order $4$ which respects the complex structure on
$\sM$.  This map lifts to an automorphism of 
 the universal cover $\tsM$ which leaves invariant a ``base point region", 
 the fundamental polycylinder defined in \S2. 
We show  in \S8 that some of the punctures in the
$c$-plane of $\sM$ are removable singularities,
so that there is a further extension of this analytic continuation
 to a multi-valued function
 on the universal cover  $\tsM^{\#}$ of  the 
 extended manifold
\beql{106aa}
\sM^{\#} := \{ (s, a, c) \in \CC \times (\CC
\smallsetminus \ZZ) \times (\CC \smallsetminus \ZZ_{\le 0})\}.
\eeq
This  manifold is obtained from $\sM$ 
by filling in the points where $c=n \ge 1$ is a positive integer.
However $\sM^{\#}$ is {\em not} preserved by the 
automorphism $\phi_{\bR}$, and the four-term functional equations cannot be extended to
be defined at all points of  its universal cover $\tsM^{\#}$. 
The universal cover $\tsM^{\#}$ appears to be a maximal
domain of holomorphy for the analytic continuation in three variables.
It is certainly maximal except for  the possible addition of  lower-dimensional  strata, of
real codimension at least two, see \S9 for further remarks.

Second, we observe that the function $\zeta(s, a,c)$ satisfies two differential-difference equations,
namely
\beql{106a}
\left( \frac{1}{2 \pi i} \frac{\partial}{\partial a} + c \right) \zeta(s, a, c) = \zeta(s-1, a ,c),
\eeq
and
\beql{106b}
 \frac{\partial}{\partial c}  \zeta(s, a, c) = -s \zeta(s+1, a ,c).
\eeq
This has the important consequence that it
satisfies a
linear partial differential equation
namely
\beql{107c}
\left(\frac{1}{2\pi i} \frac{\partial}{\partial a}\frac{\partial}{\partial c} + c \frac{\partial }{\partial c}\right)
\zeta(s, a, c) = -s \zeta(s, a, c).
\eeq
These three operators all act equivariantly with respect to the covering
map from the universal cover $\tilde{\sM} \to \sM$. It follows that 
all the monodromy functions satisfy these differential-difference
equations  and the linear partial differential equation \eqn{107c}.
Here  the differential operators $\frac{\partial}{\partial a},\frac{\partial }{\partial c}$
are viewed as acting on a complex domain,
in  the sense of Hille \cite{Hi76}.

Third, we study properties
of the monodromy functions. We
obtain the analytic continuation of the functional equations given in part I, and show
that these imply a system of linear dependencies among the monodromy functions.
We also study the vector space $\sV_s$ spanned by monodromy
functions, with the parameter $s$ held fixed.
In this case  there occur further
linear dependencies among monodromy functions, for
integer values of   $s$.
We note that at non-positive integers $s= -m \le 0$  {\em all
monodromy functions for $\zeta(s, a, c)$ vanish identically.} It follows that  the value of
the function $\zeta(s,a,c)$  is then well-defined  on the manifold $\sM$ rather than being defined
only on  a covering manifold $\tsM$.
This gives a property characterizing the these values of $s$ as ``special values".

A logical continuation of this work is to study
the change of variable $z := e^{2 \pi i a}$, which gives the
 {\em Lerch transcendent} $\Phi(s,z,c)$, given by
\beql{229}
\Phi(s,z,c) := \sum_{n=0}^\infty \frac{z^n}{ (n+c)^{s}}.
\eeq
This is done in part III, where 
we  analytically continue the function $\Phi(s, z, c)$ to a  (nearly) maximal
domain of holomorphy, and determine its monodromy functions.
The Lerch transcendent
 is closely related to the polylogarithm $Li_m(z)= \sum_{n=1}^{\infty} \frac{z^n}{n^m}$
under the specialization $c=1$ and $s=m$ is a positive integer.
We use this  to  obtain results on generalized polylogarithms and their
monodromy.

There has been very extensive prior work on analytic continuation
of the Lerch zeta function, some of it  given in
terms of the Lerch transcendent, described further in part III.
In particular  there are many individual results on analytic continuations in various
subsets of the variables, including continuations on various
 ``singular strata" above, for example that of the Hurwitz zeta function.
Most relevant to this work,  in 2000 Kanemitsu, Katsurada and
Yoshimoto \cite[Theorem 1 and Theorem 4]{KKY00} 
 obtained  an analytic continuation of the
Lerch zeta function and
Lerch transcendent in three complex variables, to a single-valued
function on various large (non-maximal)  simply-connected domains in $\CC^3$.
In \S2 we compare methods; here our objective is to  go further and determine completely the
multivalued nature of the analytic continuation. 
 We  recently discovered that in 1906 E. W. Barnes \cite{Ba06} discussed analytic
continuation of the Lerch transcendent $\Phi(s, z, c)$, and noted features of
its multivalued nature. His work appears to  give another 
approach to analytic continuation of these functions in three variables.\\

For general treatments of the Lerch zeta function we refer to the books of
 Lauren\v{c}ikas and Garunk\v{s}tis \cite{LG02},
Srivastava and Choi \cite[Chap. 2]{SC01} and Kanemitsu and Terada \cite[Chaps. 3-5]{KT07}.\\

\paragraph{\bf Acknowledgments.}
We thank  K. Prasanna and P. Sarnak
for helpful comments on this paper and P. Scott for a useful
discussion. This work was initiated while the first author worked
 at AT\&T Labs-Research and the second author consulted there. We
thank AT\&T for its support.

%
\section{Summary of results}
\setcounter{equation}{0}

The main results of this paper concern the multivalued analytic
continuation of the Lerch zeta function, summarized in Theorems \ref{th11}-
\ref{th23} below.
A more detailed  description of the multivaluedness is
given in theorems in \S4.

To effect the analytic continuation we use four
ingredients:  (i) the series representation \eqn{101}, (ii) the integral representation
\eqn{107} below, (iii) the functional equations, both
the  four-term symmetric functional equations of part I, and the
three-term asymmetric functional equation (Lerch transformation formula),
analytically continued to a suitable domains, and (iv)
the difference-differential equation \eqn{106a}. 
In comparison, the analytic continuation of Kanemitsu, Katsurada and Yoshimoto \cite{KKY00}
made use of (i), (ii) above for the Lerch transcendent; they also used the
Lerch transformation formula  and various
other expansions not considered here. 

    The series representation
\eqn{101}  defines $\zeta (s,a,c)$ as a single-valued  analytic function
of three variables  on the
polycylinder
\beql{105}
\sU  := \{ s:  s \in \CC \} \times
\{ a: \Im (a) > 0 \} \times \{c: \Re (c) > 0 \}  \subset \CC^3~.
\eeq

\noindent We also use the following  integral representation of the
Lerch zeta function:
 \beql{107}
 \zeta (s,a,c) := \frac{1}{\Gamma (s)}
\int_0^\infty \frac{t^{s-1} e^{-ct}}{1- e^{2 \pi ia} e^{-t}} dt ~.
\eeq This defines a single-valued analytic function on the
polycylinder 
\beql{108} 
\sU^+ := \{s: \Re (s) > 0 \} \times \{a: \Im(a) > 0 \} \times \{c: \Re (c) > 0 \} ~, 
\eeq 
which is smaller
than the region \eqn{105}. However this integral can be used to give
an analytic continuation to a region in the $a$-variable allowing
negative imaginary part, if the real part of $a$ is suitably
restricted. Erdelyi \cite[p. 27]{Er53} notes that in terms of the
variable $z= e^{2 \pi ia}$ the integral  above converges on the
region $z \in \CC \smallsetminus \RR_{\ge 1}$. In particular the
right side of \eqn{107} defines $\zeta (s,a,c)$ as a single-valued
analytic function on the {\em fundamental polycylinder}
 \beql{109}
 \Omega := \{ s:  0 < \Re (s) < 1  \} \times \{a : 0 < \Re (a) < 1\} \times
\{c: 0 < \Re (c) < 1 \} ~.
\eeq

      The Lerch zeta function satisfies a two four-term functional
equations, in the real variables $0< a<1, 0 < c<1$ given in part I
\cite[Theorem 2.1]{LL1}. These functional equations were obtained by
Weil \cite[p. 57]{We76}.
 Our first result, Theorem ~\ref{th11} below,
gives analytically continued versions of these functional equations
valid on  the fundamental polycylinder, and permits a
 a further  analytic continuation  in the $s$-variable to
a single-valued function on the larger region
\beql{109b}
\tilde{ \Omega} := \{ s: s \in \CC\} \times \{a : 0 < \Re (a) < 1\}
\times \{c: 0 < \Re (c) < 1 \},
\eeq
which we call the {\em extended fundamental polycylinder}.
We recall the symmetrized Lerch zeta functions 
\beql{109c}
L^{\pm}(s, a,c) := \zeta(s, a,c ) \pm e^{2 \pi i a} \zeta(s, 1-a, 1-c),
\eeq
and we let $\hat{L}^{\pm}(s, a,c)$ denote the same functions
times an appropriate archimedean factor specified in the
following result. 

\begin{theorem}\label{th11}
{\rm (Lerch Functional Equations on Polycylinders).}

(1) On the fundamental polycylinder $\Om$, which requires $0 < \Re(s) < 1$,
  the function
\beql{110}
\hat{L}^{+}(s,a,c) := \pi^{- \frac {s}{2}} \Gamma (\frac {s}{2}) \left(
\zeta (s,a,c) + e^{-2 \pi ia} \zeta (s, 1-a, 1-c) \right)~
\eeq
is holomorphic in all variables and satisfies the functional equation
\beql{110b}
 \hat{L}^{+}(s,a,c) =  e^{-2 \pi ia c} \hat{L}^{+}(1-s, 1-c, a).
 \eeq
In addition the  function
\beql{111}
L^{-}(s,a,c) := \pi^{- \frac {s+1}{2}} \Gamma (\frac {s+1}{2})
\left( \zeta (s,a,c) - e^{-2 \pi ia} \zeta (s, 1-a, 1-c)
\right)~
\eeq
is holomorphic in all variables and satisfies the functional equation
\beql{111b}
\hat{L}^{-}(s,a,c) = i e^{- 2 \pi ia c} \hat{L}^{-}(1-s, 1-c, a).
\eeq

(2) The function $\zeta (s,a,c)$ analytically continues to a holomorphic
function of three variables on the extended fundamental polycylinder
$\tilde{\Omega}$, which allows $s \in \CC$.
Both functions $L^{\pm}(s, a,c) $  analytically continue to $\tilde{\Omega}$
as holomorphic functions of three variables and the
 functional equations \eqn{111} and \eqn{111b} hold on $\tilde{\Omega}.$
\end{theorem}

Theorem ~\ref{th11} is proved in \S3.
The three-term functional equation of the Lerch zeta function
(Lerch's transformation formula) also holds on
$\tilde{\Omega}$  and is derived as Corollary 3.1.

The extended fundamental polycylinder is a large
convex domain in $\CC^3$ on which the Lerch zeta function is single-valued.
The Dirichlet series \eqn{101} representation
defines another  convex domain
\eqn{105} on which it is single-valued. The union of these two domains
is  simply connected but is not convex.

The transformation $(s,a,c) \mapsto (1-s, 1-c,a)$
 is an automorphism of the extended fundamental
polycylinder of period four. Iterating it yields
\beql{114}
\hat{L}^{\pm} (s,a,c) = (-1)^{k} e^{- 2 \pi i a} \hat{L}^{\pm} (s,
1-a, 1-c),~ for~ k~=~0~,1~ \eeq and \beql{115} \hat{L}^{\pm}
(s,a,c) = (- i)^k e^{-2 \pi i a c + 2 \pi i c} \hat{L}^{\pm} (1-s,
c, 1-a) ~,~~ for~ k~=~0,~1.
 \eeq
Here the notation $L^{+}$ corresponds to $k=0$ and $L^{-}$ corresponds to
$k=1$, i.e. $\pm= (-1)^k$.

Our next two results concern the multivalued analytic
continuation of the Lerch zeta function. Let
$\sM = \CC \times (\CC \smallsetminus \ZZ ) \times ( \CC \smallsetminus \ZZ )$,
and
let $\tilde{\sM}$ denote its universal cover, identified with
homotopy classes of curves starting from the base point
$(\frac{1}{2}, \frac{1}{2} , \frac{1}{2})$. This basepoint is in the
extended fundamental polycylinder, which
is canonically embedded  as a subset
of the universal cover $\tilde{\sM}$. We let $Z ([\gamma])$ denote
the analytic continuation of $\zeta (s,a,c)$ along a path  $\gamma
: [0, 1] \rightarrow  \sM$ in $(s,a,c)$-space which has base point
$\gamma (0) = ( \frac{1}{2}, \frac{1}{2} , \frac{1}{2})$, and
whose other endpoint is at $\gamma(1)=(s,a,c)$.
We let $[\gamma]$ denote the
homotopy class of the path $\gamma$ with fixed endpoints on the
manifold $\sM$. We will establish that the value $Z ([\gamma])$
depends only on this homotopy class. We view $Z([\gamma])$ as an
analytic function element when we allow the endpoint $\gamma (1)$
to vary, and in this case we  write $Z (s,a,c, [\gamma])$
to explicitly indicate the dependence on the endpoint. The
function, $Z (s,a,c, [\gamma])$  represents a branch of the Lerch
zeta function lying  above $(s, a, c)$  (Here we are considering a
covering manifold of a $3$-dimensional complex manifold, so
these branches are locally $3$-dimensional complex manifolds.)
In \S4 we prove the following result.
\begin{theorem}\label{th13}
{\em (Lerch Analytic Continuation)}
The Lerch zeta function $\zeta (s,a,c)$ analytically continues to
a single-valued holomorphic function $Z(s,a,c, [\gamma] )$ on the
universal cover $\tsM$ of the manifold
$$
\sM =\{ (s, a,c)  \in  \CC \times ( \CC \smallsetminus \ZZ ) \times
( \CC \smallsetminus \ZZ )\},
$$
 in which $[\gamma ] \in \pi_1 ( \sM, \bx_0)$
for $\bx_0 = ( \fot , \fot, \fot )$ and
$ \pi ([\gamma]) = (s,a,c)$ denotes the image of the covering map $\pi:
\tsM \to \sM$. Furthermore the
 function $Z(s,a,c, [\gamma] )$
is single-valued on the maximal abelian covering manifold
$\widetilde{\sM}^{ab}$ of $\sM$.
\end{theorem}

Given a multivalued function $F$ on $\sM$ and two loops $\tau_1$ and
$\tau_2$ in $\sM$ based at the point $\bx_0 = (\fot, \fot, \fot )$,
we define the {\em monodromy function}  $M(F,[\tau_1], [\tau_2])$ to
be the multivalued function whose value  on a path $\gamma$ having
$\gamma(0)=\bx_0$ is given by
 \beql{119}
  M(F,[\tau_1], [\tau_2])([\gamma]) := F( [\tau_2 \gamma ]) - F([\tau_1 \gamma]).
\eeq 
The name  ``monodromy function"
reflects that  the standard definition of  {\em monodromy}   
is  the value of this function $M(F, [\tau_1], [\tau_2])$
at the base point $\bx_0$. 

To determine all monodromy functions it suffices to 
compute them in the special case that $[\tau_1]$ is
the identity, in which case we use the simplified notation
$M_{[\tau_2]}(F):= M(F, [\tau_2], [Id])$.
In \S4 we  derive explicit formulas for the monodromy functions
of the function $Z$,
for a set of generators of $\pi_1(\sM, {\bf x_0})$, see Theorem~\ref{th31}.
We show that the monodromy functions
vanish identically on the derived subgroup (commutator subgroup)
\beql{Z120}
 (\pi_1 ( \sM , x_0))^{'} := [\pi_1 (\sM, \bx_0 ) : \pi_1( \sM, \bx_0 ) ] ~,
\eeq see Theorem~\ref{th32}. This implies that 
$Z(s,a,c , [\gamma ])$ 
defines a single-valued function on the maximal abelian covering
manifold $\widetilde{\sM}^{ab}$ of $\sM$, as stated above.

The monodromy computations reveal that all monodromy functions corresponding
to loops around the integer points $c= n \ge 1$ vanish. 
These points have removable singularities, and we obtain the
following extended analytic continuation, whose proof is deferred to \S8.

\begin{theorem}\label{th23}
{\em (Lerch Extended Analytic Continuation)}
The Lerch zeta function $\zeta (s,a,c)$ analytically continues to a
single-valued holomorphic function $Z(s,a,c, [\gamma] )$ on the
universal cover $\tilde{\sM}^{\#}$ of the manifold
$$
\sM^{\#} := \{ (s, a,c) \in \CC \times ( \CC \smallsetminus \ZZ ) \times
( \CC \smallsetminus \ZZ_{\le 0} ) \},
$$
 in which $[\gamma ] \in \pi_1 ( \sM^{\#}, \bx_0)$
for $\bx_0 = ( \fot , \fot, \fot )$ and 
$ \pi^{\#} ([\gamma]) =(s,a,c)$ denotes the image of the covering map 
$\pi^{\#}:  \tsM^{\#} \to \sM^{\#}$. Furthermore the
 function $Z(s,a,c, [\gamma] )$
is single-valued on the maximal abelian covering manifold
$(\tilde{\sM^{\#}})^{ab}$ of $\sM^{\#}$.
\end{theorem}

In \S5- \S7 we determine various  properties of the Lerch monodromy
functions, associated to the manifold $\sM$. First, in  \S5 we
observe that the Lerch zeta function satisfies the two
differential-difference equations \eqn{106a} and \eqn{106b}, under
analytic continuation.
From this we deduce that it  also satisfies
a linear partial
differential equation  \eqn{107} in the $(s, a, c) $-variables.
The  differential-difference  equations and the partial differential
equation are equivariant under the action of $\pi_1(\sM, \bx_0)$, so
we deduce  that the monodromy functions satisfy the same
difference-differential equations and partial differential equation
(Theorem~\ref{th42}). Second, in \S6 we specify linear dependencies
among the Lerch monodromy functions that follow from the functional
equations (Theorem~\ref{thK41}).  Third, in \S7 we study the vector
space $\sV_s$ generated by the Lerch monodromy functions when the
variable $s$ is fixed. Here each monodromy function $M_{[\tau]}
([\gamma])$ corresponds to some  element of $\sV$. Note that the
associated {\em monodromy representation} describing the
multivaluedness is  an induced action
$$\rho : \pi_1(\sM, {\bx_0}) \rightarrow GL(\sV_s),$$
where $\sV_s$ is a
(generally infinite-dimensional) complex vector space spanned by the germs of
function elements of $Z$ above ${\bx_0}$.
We obtain a family of representations
$\rho_s:  \pi_1(\sM, {\bx_0}) \rightarrow GL(\sV_s).$
The size of the vector space $\sV_s$ can change radically
as $s$ varies, it degenerates on certain lower dimensional strata. 
In particular we observe that for $s \in \ZZ$
there are {\em extra linear dependencies}, not implied by the functional equations, satisfied
by the monodromy functions.  We show  that the monodromy functions  vanish identically when
$s= -m \le 0 $ is a nonpositive  integer, so that
$\zeta(-m, a, c)$ is a single-valued function of $(a, c) \in (\CC \smallsetminus
\ZZ) \times (\CC \smallsetminus \ZZ)$, cf. Theorem~\ref{th41}.
We explicitly determine a basis of the vector space $\sV_s$ for all complex $s$.
For $s$ not a nonpositive integer, the vector space $\sV_s$  is infinite-dimensional.

In section \S8,  we give the  analytic continuation of
$\zeta(s, a, c)$ to
the extended manifold $\sM^{\#}$,  and prove Theorem ~\ref{th23}.\\

\paragraph{\bf Notation.}
The hat notation
$\hat{L}$ always denotes a ``completed''  function multiplied
by an appropriate archimedean Euler
factor. 
In the following  $\Re(s)$ and $\Im(s)$ denote
the real and imaginary parts of a complex variable $s$, respectively.

%
\section{Functional equations}
\setcounter{equation}{0}

 As  a preliminary step,  we recall the two functional equations
 derived in part I, originally due to A. Weil \cite{We76},
  valid for real $(a,c)$ with $0<a<1,  0<c<1$.
 We prove Theorem~\ref{th11}, analytically continuing them to the 
fundamental polycylinder $\Omega$, noting
 that $\Omega$ is invariant under the symmetries of the functional
equations.

For all real $a$ and for real $c$ with $0 < c < 1$, set
\begin{eqnarray}\label{N201}
L^{+}(s,a,c) & := & \zeta (s,a,c) + e^{- 2 \pi ia} \zeta (s, 1-a, 1-c) \nonumber \\
& = & \sum_{n= - \infty}^\infty e^{2 \pi ina} | n+c|^{-s} ~.
\end{eqnarray}
We also set
\begin{eqnarray}\label{N202}
L^{-} (s,a,c) & = & \zeta (s,a,c) - e^{- 2 \pi i a} \zeta (s, 1-a, 1-c) \nonumber \\
& = & \sum_{n=- \infty}^\infty e^{2 \pi ina} {\rm sgn} (n+c) |n+c|^{-s} ~,
\end{eqnarray}
in which
\beql{N203}
{\rm sgn} (x) = \left\{
\begin{array}{rll}
-1 & {\rm if} & x < 0 ~, \\
0 & {\rm if} & x =0 ~, \\
1 & {\rm if} & x > 0 ~.
\end{array}
\right.
\eeq
If $a \not\in \ZZ$  these series converge absolutely and define analytic functions
of $s$ for $\Re (s) > 1$; they converge conditionally (summing from $-N$ to $N$ and
letting $N \to \infty$) and define analytic functions of $s$ for $\Re(s)>0$.
If $a \in \ZZ$ they converge absolutely and define analytic functions of  $s$ for $\Re (s) > 1$.

\begin{lemma}\label{leN201} {\em (Lerch Functional Equations)}
For real $a,c$ with $0 < a < 1$ and $0 < c < 1$,
and with $\pm= (-1)^k$  the functions
\beql{N204}
\hat{L}^{\pm}(s,a,c) :=
\pi^{- \frac{s+k}{2}} \Gamma( \frac{s+k}{2}) L^{\pm}(s, a, c)
\quad\mbox{for}\quad k=0,1~,
\eeq
analytically continue to $s \in \CC$ as holomorphic functions of $s$,
and satisfy the functional equation
\beql{N205}
\hat{L}^{\pm}(s,a,c) = i^k e^{- 2\pi iac} \hat{L}^{\pm}(1-s, 1-c,a) \quad\mbox{for} \quad k =0,1 ~.
\eeq
\end{lemma}

\begin{proof}
This follows from part I, see \cite[Theorem 1.1]{LL1}. This result was originally
obtained by Weil~\cite[eqn. (15) on page 57]{We76}. Weil's
notation $(a,x,y,s)$ corresponds to our notation $(k,c,-a,2s)$.
His functions $S_a (x,y,s)$ defined by \cite[eqn. (10) on p. 56]{We76}
are related to our functions by
$S_0 (x,y,s) = L^{+} (2s,-y, x)$ and
$S_1 (x,y,s) = L^{-} (2s-1 , -y, x)$.
Some terms involve $-a$ where we have $1-a$ above, but Weil's
functions are defined to be periodic (mod $1$) in the $a$-variable.
\end{proof}

The  proof in part I of the functional equations above apply more
generally to all real $a$, including the cases $a=0$ and $1$  with
$0 < c \leq 1$. That is, one can define $\hat{L}^{\pm} (s,a,c)$
for $0 < c \leq 1$, with $a = 0$ or $1$, meromorphically continue
it to $s \in \CC$, and obtain an appropriate  functional equation.
In these cases $\hat{L}^{+} (s,a,c)$ has a simple pole at $s=1$
and no other singularities while $\hat{L}^{-}(s,a,c)$ is
holomorphic for $s \in \CC$. The case $\hat{L}^{+} (s,0,1)$
corresponds to the Riemann zeta function, and $\hat{L}^{+}(s, 0,
c)$ for $0<c<1$ to the Hurwitz zeta function. Part I observed that
an extension of $\hat{L}^{\pm} (s,a,c)$  to $(a, c) \in \RR \times
\RR$ existed that preserved the functional equation and gave
meromorphic functions in $s$ for each fixed $(a,c)$, at the cost that this function 
is discontinuous in $a$ and $c$ at integer values, for certain
ranges of $s$. This provides evidence that there cannot exist a
general analytic continuation in all three complex variables that
includes the points $a \in \ZZ$.\\

\paragraph{\em Proof of Theorem \ref{th11}.}
(i). The functional equations consist of four terms expressed in terms of
the Lerch zeta function by \eqn{110} and \eqn{111}, and each of the
four terms is a holomorphic function of three
complex variables in the fundamental polycylinder
$$\Omega = \{s : 0 < \Re (s) < 1 \} \times
\{a :0 < \Re (a) < 1 \} \times \{c : 0 < \Re (c) < 1 \} ~.
$$
These functional equations then hold by analytic continuation in the
three variables in the entire domain $\Omega$. Indeed treat $s$ and
$c$ as fixed, with $0 < c< 1$ and arbitrary $\{ s : 0 < \Re (s) < 1
\}$, and vary $a$. Since the relation holds for the interval $0 < a
< 1$, it analytically continues to all values of $a$ in the region
$\{ a :  0 < \Re (a) < 1 \}$. Once this is done, we now can choose
an arbitrary point $(s,a)$ in $\{s : 0 < \Re (s) < 1 \} \times\{a :
0 < \Re (a) < 1 \}$ and vary $c$. Since it holds for real $c$, $0 <
c < 1$, it analytically continues to $\{c : 0 < \Re (c) < 1 \}$ and
$\Omega$ is covered.

(ii). The integral \eqn{107} defines $\zeta (s,a,c)$ as an
analytic function of three variables in the domain
\beql{223}
\Omega^{+} :=
\{ s : \Re (s) > 0 \} \times \{a : 0 < \Re (a) < 1 \} \times
\{ c : 0 < \Re (c) < 1 \} ~.
\eeq
We  can suitably combine the functional equations for
$\hat{L}^{+} (s,a,c)$ and $\hat{L}^{-}(s,a,c)$ to eliminate
one of the four zeta function
terms that occur, to recover the Lerch transformation formula in the
slightly modified form
\beql{224}
\zeta (1-s,a,c) :=
(2\pi )^{-s} \Gamma (s)
\left[ e^{\frac{1}{2} \pi is - 2 \pi iac} \zeta (s, 1-c,a)+
e^{- \frac{1}{2} \pi is - 2 \pi iac + 2 \pi ic} \zeta (s,c,1-a) \right] ~.
\eeq
This formula is valid for $a$, $c$ real with $0 < a < 1$, $0 < c < 1$ and
for $s$ in the strip $\{s : 0 < \Re (s) <1 \}$, and is proved in
Theorem 5.4 of part I.
(The identities
$\Gamma (\frac{s}{2} ) \Gamma ( \frac{s+1}{2} ) = \sqrt{2 \pi} 2^{\frac{1}{2} -s} \Gamma (s)$
and $\Gamma (s) \Gamma(1-s)=\frac {\pi}{sin~\pi s}$
are used in the derivation of \eqn{224}.)
For $\Re (s) > 1$ both terms on the right side are jointly
analytic functions of all three variables in the region
$\Omega^{+}$, and taking this as a definition of the left side effects the
analytic continuation of $\zeta(s, a, c)$  to $\Re (s) < 0$ in all three variables.
The  functional equations of Theorem \ref{th11} are then inherited
by analytic continuation in the $s$-variable,
when the other two variables are held constant.
The formula \eqn{224}
expresses $\zeta (1-s,a,c)$ as a holomorphic function for
$\Re (s) > 0$, because all terms on the right are holomorphic there.
Thus $\zeta (s,a,c)$ is holomorphic in the region
$$
\Omega^- := \{s : \Re (s) < 1 \} \times \{ a : 0 < \Re (a) < 1 \}
\times \{c : 0 < \Re (c) < 1 \}$$ which with $\Omega^{+}$ covers the
extended fundamental polycylinder  $\tilde{\Omega}=\{ s : s \in \CC
\} \times \{a :  0 < \Re (a) < 1 \} \times \{c:  0 < \Re (c) < 1 \}
$. It remains to  establish the holomorphicity in $s$ of the
functions $L^{\pm}(s, a,c)$, for $s \in \CC$.  Since $\zeta(s, a,
c)$ is holomorphic in all variables, we must establish there are no
poles produced by their Gamma function factors at negative integers.
This now  follows from the functional equations taking $s$ to $1-s$,
which equate them to values at positive integers where the
holomorphicity follows from that of $\zeta(s, a,c)$ and $\zeta(s,
1-a, 1-c)$.
$~~$$~~$$~~$ $~~$$~~$$~~$ $~~$$~~$$~~$ $~~$$~~$$~~$ 
$~~$$~~$$~~$ $~~$$~~$$~~$ $\Box$ \\

The proof  above also
establishes the following variant of 
the Lerch transformation formula, valid on the extended
fundamental polycylinder \eqn{109b}.

\begin{theorem}\label{cor31}
{\em (Extended Lerch Transformation Formula)}
 For all $(s, a, c)$ in the extended fundamental polycylinder $\tilde{\Omega}$
 there holds
\beql{224a}
\zeta (1-s,a,c) =
(2\pi )^{-s} \Gamma (s)
\left[ e^{\frac{\pi is}{2}}e^{- 2 \pi iac} \zeta (s, 1-c,a) +
e^{- \frac{\pi is}{2}}e^{- 2 \pi i c(1-a)} \zeta (s,c,1-a) \right] ~.
\eeq
\end{theorem}

\begin{proof}  This holds for $(s, a,c) \in \tilde{\Omega}$ using
\eqn{224} .
\end{proof}

 The extended Lerch transformation formula above differs from the original
Lerch transformation formula \eqn{102} in having
the term $\zeta (s, 1-c,a)$
in place of $\zeta (s, -c, a)$.
This formula still agrees with  \eqn{102}
because the function $\zeta (s,a,c)$ in Lerch's original
paper is periodic (mod one) in the $a$-variable for
$\Im(a) > 0$, and then by extension for parts of the real axis, so that
\beql{225}
\zeta (s, -c,a) = \zeta (s, 1-c,a)
\eeq
holds in a suitable domain. The extended Lerch transformation formula \eqn{224a}
keeps all Lerch zeta function terms inside
the extended fundamental polycylinder.
Under analytic continuation the equality \eqn{225}
is preserved only if the correct
branches of the function are chosen on both sides, cf. Theorem~\ref{thK41}.

%

\section{Analytic continuation and monodromy functions}
\setcounter{equation}{0}

In this section we analytically continue $\zeta (s,a,c)$ to a
multivalued function on the manifold
\beql{Z301}
\sM := \CC \times ( \CC \smallsetminus \ZZ ) \times ( \CC \smallsetminus \ZZ ) ~,
\eeq
which is coordinatized with the three complex variables
$(s,a,c)$, and prove Theorem \ref{th13}.  We fix the base point ${\bx_0}:= (\fot, \fot, \fot)$
in the fundamental polycylinder $\Omega$ and consider analytic
continuation of the Lerch zeta function along a path $\gamma :
[0,1] \to \sM$ emanating from the base point $\gamma (0) = \bx_0$
to an endpoint $\gamma (1) = (s,a,c) \in \sM$. The analytic
continuation will depend only on the homotopy class $[\gamma ]$ of the
path (where homotopies hold endpoints fixed).
\begin{defi}
{\em The {\em multi-valued function element }$Z(\gamma)$ denotes the
 analytic continuation of the Lerch zeta
function $\zeta(s, a, c)$ along the path $\gamma$, starting from
the base point $(s, a, c) = (\fot, \fot, \fot).$  We shall
sometimes denote this by
$Z( s, a, c; \gamma)$, where $(s, a, c)$ are the coordinate values at the endpoint
of $\gamma$, and view it as a single-valued holomorphic function in a small open neighborhood
of the endpoint $(s, a,c)$.}
\end{defi}
Below we will use the notations $Z([\gamma])$ and $Z([s, a, c, [\gamma])$,
anticipating  the dependence of function elements only on the homotopy class of the path;
 this latter property is only established in the proof of
Theorem~\ref{th31} below.

We may alternatively regard $Z([\gamma])$ as a holomorphic function
on the universal covering manifold $\tsM$ of $\sM$. The
universal cover $\tsM$ is the collection of homotopy classes
$[\gamma]$ of paths from $\bx_0$, and the projection 
$\pi: \tsM \to \sM$ is the endpoint 
\beql{Z302} 
\pi ( [\gamma] ) :=\gamma (1) ~.
 \eeq 
 The lifted base point  in the universal cover is
  $\tilde{\bx}_0 := [\bx_0]  \in \tsM$,
 and we let $\tilde{\bx} = [\gamma]$ 
 denote a general point  in $\tsM$.
 Since $\tsM$
  is a covering space of $\sM$, it
inherits its complex-analytic structure. Thus if we can analytically
continue $\zeta (s,a,c)$ along a path $\gamma$, it will
automatically be holomorphic in a neighborhood of its endpoint. When
we want to view $Z ([\gamma])$ as a holomorphic function element in
a neighborhood of $[\gamma]$, we write $Z(s,a,c, [\gamma ])$. As the
extended fundamental polycylinder $\tOm$ is simply connected, we
shall identify a point $\bx \in \tOm$ with a path $\gamma$ from
$\bx_0$ to $\bx$ which remains entirely in $\tOm$. In this way we
embed $\tOm$ in $\tsM$. Given one path $\gamma$ from $\bx_0$ to
$\bx$, we may describe all homotopy classes of paths from $\bx_0$ to
$\bx$ by $[ \tau \gamma]$, where $[\tau ]$ runs through all classes
of loops in the fundamental group $\pi_1 (\sM , \bx_0)$ of $\sM$
with base point $\bx_0$.

We now specify a set $\sG$ of generators of the homotopy
group  $\pi_1(\sM, \bx_0)$.
In view of \eqn{Z301}, the $\pi_1 (\sM ,\bx_0 )$
is the product of the homotopy groups of the three
manifolds $\CC$, $\sX_a = \CC \smallsetminus \ZZ$ and
$\sY_c = \CC\smallsetminus \ZZ$. The homotopy group $\pi_1 ( \CC , \frac{1}{2})$
is trivial. The group $\pi_1 ( \sX_a, \frac{1}{2} )$ is a free
group on countably many generators $\{ [ X_n ] : n \in \ZZ \}$, in
which $X_n$ is a small counterclockwise-oriented loop around the
point $a = n \in \ZZ$ in the $a$-plane, which is initially reached
by a path from the base point $a= \frac{1}{2}$ which lies in the
upper half-plane $\{ a: \Im (a) > 0 \}$ and which returns to
$\frac{1}{2}$ along the same path. The group 
$\pi_1 ( \sY_c, \frac{1}{2} )$ is defined similarly, with generators 
$ \{[Y_n]: n\in \ZZ \} $, 
where $[Y_n]$
traverses a small counterclockwise loop  around the point $c = n$ in
the $c$-plane, reached along a path in 
$\{ c: \Im(c) >0\}$. These give the  set of generators
\beql{Z303}
 \sG := \{ [X_n ] : n \in \ZZ \}
\cup \{ [Y_n ]: n \in \ZZ \}
\eeq
of  $\pi_1 ( \sM, \bx_0)$.

We determine explicit
formulas for the multivalued function elements $Z([\gamma ])$
in terms of the action of the generating
set $\sG$, given in Theorem ~\ref{th31} below. We first 
outline  the
proof approach.
We begin by  analytically continuing $\zeta (s, a, c)$ to holomorphic functions
on a collection  $\{ \sO_j \}$ of open, simply connected subdomains
$\sO_j$ of $\sM$,
each containing $\Om$ such that their union covers $\sM$.
Hence we have extended $\zeta$ to a holomorphic function $Z$
on {\em some open subset $\, \sW$ of $\tsM$} whose projection to $\sM$
is the whole space.
Next we observe that for each point $\bx$ in the intersection of two subdomains,
we have obtained two values $Z( [\gamma_1 ])$ and $Z([\gamma_2] )$,
following two paths $\gamma_1$, $\gamma_2$ from $\bx_0$ to $\bx$ in each
subdomain.
We call the difference $Z( [\gamma_2 ]) - Z([\gamma_1 ])$ the
{\em monodromy} of
$Z$ at $[\gamma_1]$ with respect to the loop
$[\tau ] \in \pi_1 (\sM , \bx_0 )$, where
$\tau = \gamma_2 \gamma_1^{-1} $ is the composition of the path $\gamma_2$ followed
$\gamma_1$ traced backwards.
We view this as a function of the path $[\gamma_1]$, called the
{\em monodromy function}
$M_{[\tau]} (Z) ([\gamma_1 ])$, and formalize this in the following
definitions.
%
\begin{defi}\label{de31}
{\rm Let $f: \tsM \to \CC$ be a continuous function on the
universal cover $\tsM$ of a manifold $\sM$ and let $[\tau ] \in\pi_1 ( \sM , \bx_0 )$
 be a homotopy class. The {\em
$[\tau]$-translated function} $\sQ_{[\tau]} (f): \tsM \to \CC$
is defined by
 \beql{Z304}
\sQ_{[\tau]} (f) ( [\gamma ] ) := f([\tau \gamma ]) ~, 
\eeq 
where
$\gamma$ is a path with basepoint $\gamma (0) = \bx_0$ and $\tau  \gamma$
 is the composed path.}
\end{defi}

Each element  $[\tau ]$ of $\pi_1 (\sM, \bx_0)$ corresponds to a
fixed-point-free homeomorphism $\ot$
of $(\tsM, \tilde{\bx}_0)$ which commutes with the projection
$\pi: \tsM \to \sM$
(see Hatcher \cite[Chap. 1.3]{Ha01}). Then 
 we have, in terms of  $\ot$, 
\beql{Z304b}
\sQ_{[\tau]}(f)(\tilde{\bx}) = f(\ot  (\tilde{\bx})).
\eeq
Each  $\sQ_{[\tau]}$ acts as a linear operator on the vector space $\rC^0 (\tsM )$ of
continuous complex-valued functions on $\tsM$, defining an  representation
of the group $\pi_1 (\sM, \bx_0)^{opp}$ having the opposite multiplication, i.e.
\beql{Z305}
\sQ_{[\tau_2]} \sQ_{[\tau_1]} = \sQ_{[\tau_1 \tau_2]} ~.
\eeq
We have $\sQ_{[\tau^{-1}]} = \sQ_{[\tau]}^{-1}$ and    the operators  
$\tilde{\sQ}_{[\tau]} :=\sQ([\tau])^{-1}$ 
 give a linear representation
of $\pi_1 (\sM, \bx_0)$ on the vector space $\rC^0 (\tsM )$.


\begin{defi}\label{de32}
{\rm
Let $f: \tsM \to \CC$ be a continuous function and let
$[\tau ] \in \pi_1 (\sM , \bx_0)$ be a homotopy class.
The {\em monodromy function} $M_{[\tau]} (f): \tsM \to \CC$
of $f$ at $[\tau ]$ is defined by
\beql{Z305A}
M_{[\tau]} (f) := ( \sQ_{[\tau]} - I ) (f) ~.
\eeq
That is, for all paths $\gamma : [0,1] \to \sM$ with base point
$\gamma (0) = \bx_0$,
\beql{Z306}
M_{[\tau]} (f) ([\gamma]) :=
f([ \tau \gamma]) - f ([\gamma ]) ~.
\eeq
}
\end{defi}

Monodromy functions obey the following
relation.


\begin{lemma}\label{le31}
For a single-valued function $f(\tilde{\bx})$ on $\tsM$ and any
$[\tau_1 ]$, $[\tau_2 ] \in \pi_1 (\sM , \bx_0)$, we have
\beql{Z307}
M_{[\tau_1 \tau_2 ]} (f) =
M_{[\tau_1]} (f) +
M_{[\tau_2]} (f) + M_{[\tau_2]} (M_{[\tau_1]} (f)) ~.
\eeq
\end{lemma}

\begin{proof}
We have
$$M_{[\tau_1 \tau_2]}(f) ([\gamma ] )= f( [ \tau_1 \tau_2\gamma \ ]) - f([\gamma] ) ~,
$$
and
$$M_{\tO} (f) ([\gamma] ) + M_{\tT} (f) ([\gamma ]) =
f([\tau_1 \gamma  ]) - f( [\gamma] ) +
f([\tau_2 \gamma  ]) - f ( [\gamma ] ), ~
$$
while
\begin{eqnarray*}
M_{[\tau_2]} (M_{[\tau_1]} (f))( [\gamma ]) & = & M_{[\tau_1]} (f) ([\tau_2 \gamma  ]) -
M_{[\tau_1]} (f) ([\gamma]) \\
& = & (f([\tau_1 \tau_2\gamma  ]) - f([\tau_2\gamma  ])) -
(f( [\tau_1 \gamma  ]) - f( [\gamma] )) ~.
\end{eqnarray*}
Combining these formulae yields \eqn{Z307}.
\end{proof}

To continue outlining the proof, we obtain a partially defined
monodromy function
$M_{[\tau]} (Z) ([\gamma ])$ defined only for certain paths $\gamma$
lying in one subdomain  $\sO_j$  and ending in the overlapped area.
However the monodromy functions $M_{[\tau]} (Z)$, each defined on
an open subset of $\tsM$, are ``simpler'' functions than $Z$, and
can themselves be {\em directly} analytically continued to $\tsM$
without knowledge of $Z$. We  are then able to analytically
continue $Z$ along all loops
 $[\tau ] \in \pi_1 ( \sM , \bx_0)$,
by writing $[\tau ]$
as a word in the generating set $\sG$, and repeatedly applying
Lemma \ref{le31} to compute $M_{[\tau]} (Z)$ and finally obtain
\beql{Z308}
Z([ \tau\gamma ]) := M_{[\tau]} (Z) ([\gamma]) +
Z([\gamma]) ~,
\eeq
for $[\gamma ]$ contained in $\sO_j$, and
$[\tau ] \in \pi_1 ( \sM , \bx_0 )$. Since the sets $\{ \sO_j \}$
cover $\sM$, this defines a single-valued holomorphic function $Z$
on $\tsM$, which agrees with $\zeta (s, a, c)$ on the fundamental
polycylinder. The resulting  function $Z$ is
 the desired function by uniqueness of
analytic continuation along paths, using 
Cauchy's theorem and the fact that $\tsM$ is simply connected.


The following result establishes analytic continuation and
explicitly gives monodromy functions associated to the generating
set $\sG$. In this result,  the {\em principal branch} of the
function $\log z$ is the branch of $\log z$ on $\CC^*$ which is real
on the positive real axis, with the branch defined by making a cut
along the negative imaginary axis, so that on the negative real axis
one has $\log (-x) := \log x + \pi i$, where $x \in \RR_{>0}$.

%
%

\begin{theorem}\label{th31}
{\em (Lerch monodromy functions)}
The Lerch zeta function $\zeta (s, a, c)$ analytically continues
to a single-valued
holomorphic function $Z$ on the universal cover $\tsM$ of $\sM$.
The monodromy functions of $Z$ on the generators of $\pi_1 (\sM , \bx_0)$
are given as follows.

(i-a)~For $n \in \ZZ$, on the extended fundamental polycylinder
$\tilde{\Omega}$ given by {\rm \eqn{109b}} we have
\beql{WL310}
M_{[X_n]} (Z) (s,a,c) = 
-\frac{(2 \pi)^s e^{\frac{\pi is}{2}}}{\Gamma (s)} (a - n)^{s-1} e^{- 2 \pi ic (a-n)},
 \eeq
  where
$ (a-n)^{s-1} := e^{(s-1) \log (a-n)}$ using the principal
branch of the logarithm. This function analytically continues to a
single-valued function on  $\tsM$.

(i-b) For $n \in \ZZ$, on any path
$\gamma$  in $\sM$ from $\bx_0$ to an endpoint $\bx_1$ lying in
the multiply connected region
\beql{WL311aa}
 \sM_s :=\{s\}\times ( \CC \smallsetminus \ZZ ) \times ( \CC \smallsetminus \ZZ)
\eeq
we have
\beql{WL311}
M_{[X_n]^{-1}} (Z) ([\gamma])=  - e^{-2 \pi is}M_{[X_n]} (Z)([\gamma])~,
 \eeq
 and
 \beql{WL312}
  M_{[X_n]^{\pm k}}(Z) ([\gamma])= \frac{e^{\pm 2 \pi isk} -1}{e^{\pm 2\pi is} -1}
M_{[X_n]^{\pm 1}} (Z) ([\gamma])\quad\mbox{for} \quad k \ge 1 ~.
\eeq

(ii-a) For $n \in \ZZ$, on the extended fundamental polycylinder $\tilde{\Om}$
there holds
\beql{WL313}
M_{[Y_n]} (Z)(s, a, c, [\gamma]) =
\left\{
\begin{array}{cl}
0 & \mbox{for}~~ n \ge 1, \\
(e^{- 2\pi is} -1) e^{- 2 \pi ina} (c-n)^{-s}   & \mbox{for}~~ n \le
0,
\end{array}
\right.
\eeq
 where $(c-n)^{-s} :=e^{-s \log (c-n)}$ using the principal branch of the logarithm. This
function analytically continues to a single-valued function on
$\tsM$.

(ii-b) For $n \in \ZZ$, on any path $\gamma$ in $\sM$ from $\bx_0$ to an endpoint
$\bx_1$ lying in $\sM_s$ we have
we have
\beql{WL314}
M_{[Y_n]^{-1}} (Z) ([\gamma])= -e^{2 \pi is} M_{[Y_n]} (Z) ([\gamma])~,
\eeq
and
\beql{WL315}
M_{[Y_n]^{\pm k}} (Z)([\gamma]) =
 \frac{e^{\mp 2 \pi i s k}-1}{e^{\mp 2 \pi is} -1} M_{[Y_n]^{\pm 1}} (Z)([\gamma])
\quad\mbox{for}\quad k \ge 1 ~.
 \eeq
\end{theorem}


\paragraph{\em Remark.} With our convention on the principal
branch of the logarithm, the formula \eqn{WL310} has the more explicit form
\begin{equation}~\label{WL310a}
M_{[X_n]} (Z) (s,a,c) =\left\{
\begin{array}{lll}
-\frac{(2 \pi)^s e^{\frac{\pi i s}{2}}}{\Gamma (s)}
e^{\pi i (s - 1)}(n - a)^{s-1} e^{- 2 \pi ic (a-n)} ~ & {\rm if} & n \geq 1,
\nonumber \\
 & & \nonumber \\
-\frac{(2 \pi)^s e^{\frac{\pi i s}{2}}}{\Gamma (s)}
(a-n)^{s-1} e^{- 2 \pi ic (a-n)} ~~~~~  & {\rm if} & n \leq 0. \nonumber \\
\end{array}
\right.
\end{equation}

\begin{proof}
 Recall that the extended fundamental polycylinder
$\tilde{\Omega}$ is embedded as a subset of $\tilde{\sM}$ with
base point $(\frac {1}{2}, \frac {1}{2}, \frac {1}{2})$
and that on this region the function $Z$
agrees with $\zeta(s, a, c)$.

We start with the integral representation
of the Lerch zeta function
\beql{Z325}
\zeta (s,a,c) := \frac{1}{\Gamma (s)} \int_0^\infty \frac{e^{-ct}}{1-e^{2 \pi ia} e^{-t}} t^{s-1} dt~.
\eeq

Let $\sA_L$ denote the $a$-plane cut along all the half-lines
 $M_k = \{ k- it : 0 \le t < \infty \}$ for $k \in \ZZ$, so that
\beql{WL317}
 \sA_L := \CC \smallsetminus \{ M_k : k \in \ZZ \} ~.
\eeq
This cut region is pictured in Figure \ref{fg31}.
\begin{figure}[htb]
\begin{center}
\input cut1.pstex_t
\end{center}
\caption{Cut domain $\sA_L$.}
\label{fg31}
\end{figure}

The integral \eqn{Z325} defines a single-valued function on the
simply-connected domain
\beql{Z327}
\sO_1 :=
\{ s: \Re (s) > 0 \} \times
\{ a: a \in \sA_L \} \times \{ c: \Re (c) > 0 \} ~.
\eeq
Furthermore we may obtain a similar result for a domain of altered shape by
considering instead
\beql{Z328}
\zeta (s,a,c) := \frac{1}{\Gamma (s)} \int_{L_{u, \epsilon}}
\frac{e^{-ct}}{1- e^{2 \pi ia} e^{-t}} t^{s-1} dt ~,
\eeq
in which $L_{u,\epsilon}$ is a deformed path from 0 to $\infty$ in the
$t$-plane on the real axis
which makes a clockwise detour around the point $u \in \RR_{>0}$
through a circle of radius $\ep$, with $0 < \ep < min (u, 1/2)$,
 as shown in Figure \ref{fg32}.
\begin{figure}[htb]
\begin{center}
\input contour.pstex_t
\end{center}
\caption{Contour $L_{u,\ep}$.}
\label{fg32}
\end{figure}

The integral \eqn{Z328} defines a single-valued holomorphic function
on the simply-connected domain
$$
\sO_1 (u, \ep ) := \{
s: \Re (s) > 0 \} \times \sA_{L_{u, \ep}} \times
\{ c: \Re (c) > 0 \}
$$
in which $\sA_{L_{u,\ep}}$ denotes the $a$-plane cut along all
the deformed half-lines
$$
L_k^- (u,\ep ) :=
\left\{
a = k-
\frac{i}{2 \pi} t:
t \in L_{u,\ep} \right\} \quad\mbox{for}\quad k \in \ZZ ~,
$$
so that
\beql{Z329}
 \sA_{L_{u, \ep}} := \CC \smallsetminus \{L_k^- (u,\ep ) : k \in \ZZ \} ~.
\eeq
As $u$ and $\ep$ vary over
$0 < u < \infty$ and $0 < \ep < min (u, 1/2)$, the domains $\sO_1
(u,\ep )$ cover the entire domain
\beql{Z330}
\sO_1^\infty := \{s: \Re (s) > 0 \} \times \{ a: a \in \CC \smallsetminus \ZZ \}
\times \{ c: \Re (c) > 0 \} ~.
 \eeq
  Moreover, for $u > 0$ and $0 < \ep < min (u, 1/2)$, the open half-disk
 \beql{Z331}
 D_{u,\ep} (n) := \left\{ a:  \left| a- \left( n- \frac{iu}{2 \pi} \right)\right|
 < \ep \quad\mbox{and}\quad \Re (a) > n \right\}
\eeq
lies
in the intersection of the domains $\sA_{L_{u,\ep}}$ and $\sA_L$,
but in homotopically different components, see Figure \ref{fg33}.
We use this fact to compute $M_{[X_n]} (Z)$, where $Z$ denotes the
analytic continuation of the function $\zeta (s,a,c)$ from the
base point $x_0$. We can find two paths $\gamma_1$, $\gamma_2$
with base point $x_0 = ( \frac{1}{2}, \frac{1}{2}, \frac{1}{2} )$
which hold $c=s=1/2$ fixed and vary $a$, such that
$\gamma_1 (1) =\gamma_2 (1) = ( \frac{1}{2} , a , \frac{1}{2} )$
with $a \in D_{u,\ep} (n)$, while  $\gamma_1 ([0,1]) \subseteq \sA_L$
and $\gamma_2 ([0,1]) \subseteq \sA_{L_{u,\ep}}$.
Then $\gamma_2 \gamma_1^{-1} $ is homotopic to $X_n$ in $\sM$. Hence
 \beql{Z332}
Z([\gamma_2]) - Z([\gamma_1]) = Z([(\gamma_2 \gamma_1^{-1})\gamma_1])
 - Z([\gamma_1 ]) = M_{[X_n]} (Z) ([\gamma_1])~.
\eeq
The contours $\gamma_1$ and $\gamma_2$ are pictured in Figure
\ref{fg33}.
\begin{figure}[htb]
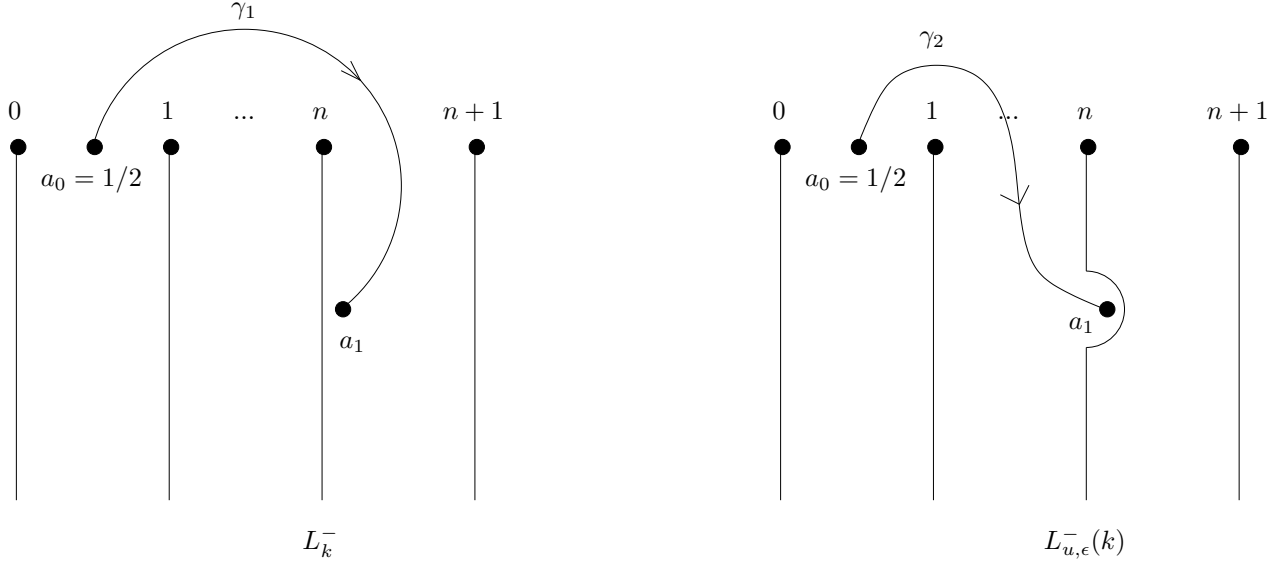

\begin{center}
\input 8lines.pstex_t
\end{center}
\caption{Contours for $M_{[X_n]} (Z)$.}
\label{fg33}
\end{figure}

Now we have
\begin{eqnarray*}
Z([\gamma_2]) & = & \frac{1}{\Gamma (s)} \int_{L_{u,\ep}}
\frac{e^{-ct}}{1-e^{2 \pi ia} e^{-t}} t^{s-1} dt \\
& = & \frac{1}{\Gamma (s)}
\int_0^\infty \frac{e^{-ct}}{1-e^{2 \pi ia} e^{-t}} t^{s-1} dt +
\frac{1}{\Gamma (s)} \oint_{C_{u,\ep}}
\frac{e^{-ct}}{1-e^{2 \pi ia} e^{-t}} t^{s-1} dt ~,
\end{eqnarray*}
in which $C_{u, \ep}$ is a closed clockwise oriented semicircular
contour centered at $u$ with radius $\ep$ together with the
negatively oriented line segment $[u-\ep, u+\ep ]$. The residue
theorem gives 
\beql{Z333} 
Z([\gamma_2 ]) = Z([\gamma_1]) - \frac{2\pi i}{\Gamma(s)} 
Res_{t= 2\pi i(a-n) } \left( \frac{e^{-ct}}{1-e^{2 \pi ia} e^{-t}} t^{s-1} \right) ~,
 \eeq
  in which 
$a= ( n-\frac{iu}{2 \pi} ) + (x+ i y )$ with $x > 0$ and $x^2 + y^2 <\ep^2$.
 We deduce that 
\beql{Z334} 
M_{[X_n]} (Z) ([\gamma_1 ]) = -
\frac{(2 \pi i)^s}{\Gamma (s)} (a-n)^{s-1} e^{- 2 \pi ic (a-n)}
\quad\mbox{for}\quad n \in \ZZ  
\eeq 
holds near the endpoint $(s,a,c)$ of
$\gamma_1$. Here 
$$(2 \pi i )^s := \exp (s(\log 2 \pi i)) = (2 \pi)^s e^{\frac{\pi is}{2}},$$
where 
$\log 2 \pi i := \log 2\pi + \frac{i\pi}{2}$ 
is the value on the principal branch of $\log s$, and
\newline
$(a-n)^{s-1} := \exp ((s-1) \log (a-n))$ with $\log (a-n)$
determined by the principal branch of the logarithm, observing that
$\Re (a-n) > 0$ at the endpoint of $\gamma_1$, and $a - n$ never
crosses the  negative imaginary axis when tracing the path
$\gamma_1$ backwards. The function
{\em defined} by the right side of \eqn{Z334} on $\sA_L$
analytically continues to a holomorphic function on the universal
cover $\tsM$, and its only branch points over $\sM$ are at 
$\CC \times \{ a = n\} \times ( \CC \smallsetminus \ZZ )$. In
particular when this function is analytically continued inside
$\sA_L$ back to the fundamental polycylinder along
$[\gamma_1]^{-1}$, we obtain the functions given in \eqn{WL310}.

Next, we use the series representation of the Lerch zeta function
\beql{WL326}
\zeta (s,a,c) = \sum_{n=0}^\infty e^{2 \pi ina}(n+c)^{-s} ~,
\eeq
which defines the principal branch of $\zeta(s,a,c)$ as a
 holomorphic function on the simply-connected domain
\beql{WL327} 
\sU = \{ s: s \in \CC \} \times \{ a: \Im (a)> 0 \} \times \{ c: \Re (c) > 0 \} ~. 
\eeq 
For integers $n \ge 0$, the function
 \beql{WL328}
 F_n (s,a,c) := e^{2 \pi in a} (n+c)^{-s} = e^{2 \pi ina} e^{- s \ln (n+c)}
\eeq
on  $\sU$ analytically continues to a multi-valued function on $\sM$
whose branch locus is 
$\CC \times ( \CC \smallsetminus \ZZ ) \times\{ c = -n\}$. 
Thus all the monodromy functions 
$M_{[X_l]} (F_n) =0$ ($l \in \ZZ$), while 
\beql{WL329} 
M_{[Y_l]} (F_n ) = 
\left\{
\begin{array}{cll}
(e^{- 2 \pi is} -1) F_n & \mbox{if} & l = -n ~, \\ [+.1in]
0 & \mbox{if} & l \neq -n ~,
\end{array}
\right.
\eeq
holds in the submanifold 
$\{s\} \times (\CC \smallsetminus\ZZ ) \times (\CC \smallsetminus \ZZ )$.
 By separating the terms
$\sum_{n=0}^N e^{2 \pi ina} (n+c)^{-s}$, 
we may rewrite \eqn{WL326} as
\begin{eqnarray}\label{WL330}
\zeta (s,a,c) & = & \sum_{n=0}^N e^{2 \pi ina} (n+c)^{-s} +
\sum_{n = N+1}^\infty e^{2 \pi ina} (n+c)^{-s} \nonumber \\
& = & \sum_{n=0}^N F_n (s,a,c) + e^{2 \pi i (N+1) a} \zeta (s,a,c+N+1) ~.
\end{eqnarray}
We have already shown that the last term
$\zeta (s,a,c+N+1)$ can be holomorphically continued
to simply connected regions covering the domain
$$
\sO_1^\infty (N) := \{ s: \Re (s) > 0 \} \times \{ a: a \in \CC
\smallsetminus \ZZ \} \times \{ c: \Re (c) > -N-1 , c \not\in \ZZ_{\le 0} \}.
$$
and the same holds for  each of the $F_n$. Therefore the right side of \eqn{WL330}
defines an analytic continuation  for the left side $\zeta(s,a,c)$.
 Letting $N \to \infty$, we conclude that $\zeta (s,a,c)$
can be holomorphically continued to simply connected regions
covering the domain 
\beql{WL331}
 \sM^+ := \{ s: \Re (s) > 0 \} \times \{ a: a \in \CC \smallsetminus \ZZ \} \times \{ c : c \in
\CC \smallsetminus \ZZ_{\le 0} \} ~. 
\eeq 
Finally we use the extended Lerch
transformation formula given  in Theorem \ref{cor31} 
 to extend $\zeta (s,a,c)$ to a
holomorphic function on
simply connected regions covering the domain
$$
\sM^- := \{s : \Re (s) < 1 \} \times
 \{ a: a \in \CC \smallsetminus \ZZ \} \times \{ c: c \in \CC \smallsetminus \ZZ \} ~.
 $$
The resulting continuation is holomorphic since $\Gamma (s)$ is
holomorphic in the half-plane $\Re (s) > 0$.

We next  compute the monodromy functions of $Z$ along $[Y_n]$,
e.g. $M_{[Y_n]} (Z)$.
A loop in $\sM$ based at $x_0$ around $c = 1,2,3, \ldots$ that
remains in the positive $c$-plane and along which
$a = s = \frac{1}{2}$ are held fixed is contractible in $\sM^+$,
hence we find that
\beql{WL332}
M_{[Y_n]} (Z) \equiv 0 \quad\mbox{for}\quad n =1,2,3 ,\ldots ~.
\eeq
In view of \eqn{WL330}, \eqn{WL329} and \eqn{WL332}, for integer $n \le 0$,
the monodromy $M_{[Y_n]} (Z)$ arises from the monodromy of $F_{-n}$.
In other words, on the extended
fundamental polycylinder $\tilde{\Omega}$,
$$
M_{[Y_n]} (Z) (s,a,c) = (e^{- 2 \pi is} -1 ) e^{- 2 \pi ina} (c-n)^{-s} ~~~\mbox{for}~~n =0, -1, -2, ...
$$
as given in \eqn{WL313}.

To complete the proof of Theorem~\ref{th31}, we verify 
the remaining parts of assertions (i)--(ii).
We have already verified the formulas \eqn{WL310} and \eqn{WL313}.
Now Lemma \ref{le31} applied to $\tau_2 = \tau_1^{-1}$ gives
\beql{437aa}
0= M_{[Id]} (Z) = M_{[\tau_1]} (Z) + M_{[\tau_1]^{-1}} (Z) +
M_{[\tau_1]^{-1}} (M_{[\tau_1]} (Z) ) ~.
\eeq
Choosing  $[\tau_1] = [X_n]$, we evaluate the last
term  on the right side using  \eqn{Z334}, to obtain
$$
M_{[X_n]^{-1}}( M_{[X_n]}(Z))  = 
M_{[X_n]^{-1}}( -\frac{(2 \pi i)^s}{\Gamma (s)} (a-n)^{s-1} e^{- 2 \pi ic (a-n)}) 
= (e^{-2\pi i s} -1)M_{[X_n]}(Z).
$$
We  deduce  from \eqn{437aa} that for a path $\gamma$
with  $\gamma(0) = \bx_0$ and  other endpoint in $\sM_s$ that
$$
M_{[X_n]^{-1}} (Z) ([\gamma])=  -e^{-2 \pi is} M_{[X_n]} (Z)([\gamma] )~.
$$
To verify the formulas for $M_{[X_n]^{\pm k}}(Z)$ we also use Lemma \ref{le31}.
We  proceed by induction on $k$. Treating first $k \ge 1$, the base case $k=1$
holds trivially and using the induction hypothesis we obtain
\begin{eqnarray*}
M_{[X_n]^{k+1}}(Z) &= &M_{[X_n]^k}(Z) + M_{[X_n]}(Z) + M_{[X_n]}( M_{[X_n]^k}(Z)).\\
&=&  \frac{e^{2 \pi i sk} -1} {e^{2\pi i s}-1}M_{[X_n]} (Z) + M_{[X_n]} (Z) 
+(e^{2\pi is}-1) \left( \frac{e^{2 \pi i sk} -1} {e^{2\pi i s}-1} M_{[X_n]}(Z) \right) \\
&=& 
=\frac{e^{2 \pi i s(k+1)} -1} {e^{2\pi i s}-1}M_{[X_n]}(Z),
\end{eqnarray*}
completing the induction step. The proof for $k \le -1$ is similar.
Finally we deduce the formulas for $M_{[Y_n]^{-1}} (Z)$ and $M_{[Y_n]^{\pm k}}$
by analogous arguments.
\end{proof}

The next result determines all monodromy functions for $Z$ and shows
that $Z$ is single-valued on the maximal abelian cover
$\tsM^{ab}$ of $\sM$.
%
%
\begin{theorem}\label{th32}
Suppose that $[\tau] \in \pi_1 ( \sM , \bx_0)$ satisfies
\beql{Z342} [\tau ] = [S_1]^{\ep_1} [S_2]^{\ep_2} \cdots
[S_m]^{\ep_m} ~, \eeq with each $[S_i]$ in the set of generators
$\sG$ of $\pi_1 (\sM , \bx_0)$, in which each $\ep_i$ equals one
of $\pm 1.$ Then the monodromy function $M_{[\tau]} (Z)$ for the
Lerch zeta function on the universal cover $\tsM$ satisfies
\beql{Z343}
M_{[\tau]} (Z) = \sum_{[S] \in \sG} M_{[S]^{k(S)}} (Z)
~, \eeq
 with
\beql{Z344}
k(S) := \sum_{\{ i: S_i = S\}} \ep_i~ ,\quad\mbox{for each} \quad S \in \sG ~.
\eeq
In particular,
$M_{[\tau]} (Z)$ vanishes identically for
$[\tau ]$ in the
commutator subgroup
\beql{Z345}
 \Gamma :=
 ( \pi_1 (\sM, \bx_0)) ^{'}
= [ \pi_1 (\sM , \bx_0): \pi_1 (\sM , \bx_0) ] 
 \eeq
  of $\pi_1(\sM, \bx_0)$. 
  Thus $Z$ is a single-valued function on the quotient manifold \\
$\tsM^{ab} := \Gamma \backslash \tsM $ 
with $\Gamma$ acting as
a group of homeomorphisms of $\tsM$.
\end{theorem}

\paragraph{\em Remark.}
The manifold $\tsM^{ab}$ is the maximal unramified abelian
covering manifold of $\sM$ with projection map $\pi_{ab}:
\tsM^{ab} \to \sM$ induced from $\pi: \tsM \to \sM$. We let
$\pi_{\bf D}: \tsM \to \tsM^{ab}$ denote the
projection from the universal cover to $\tilde{\sM}^{ab}$.

\begin{proof}
Recall that $\sG =\{ [X_n]: n \in \ZZ \} \cup \{ [Y_n]: n \in \ZZ \}$.
We first note that for any two distinct generators $[S_1]$, $[S_2] \in \sG$
the Lerch monodromy functions are independent in the sense that
\beql{Z346}
M_{[S_1]^{k_1}}( M_{[S_2]^{k_2}}( Z)) = 0 ~.
\eeq
This follows by a direct computation from Theorem~\ref{th31},
because $M_{[S_2]^{k_2}}( Z)$ has nontrivial monodromy only around the
generator
$[S_2]$ in $\sG$ and no other generator.
Now Lemma~\ref{le31} gives for $[S_1] \neq [S_2]$ that
\beql{Z347}
M_{[S_1]^{k_1} [S_2]^{k_2}}( Z) = M_{[S_1 ]^{k_1}}( Z)+
M_{[S_2 ]^{k_2}}(  Z) ~.
\eeq

We prove Theorem \ref{th32} by induction on $m$.
The base case $m=1$ holds because both sides of \eqn{Z343} are then identical;
\eqn{Z347} shows that the case $m =2$ also holds.
For the induction step, assume it holds for $m-1 \ge 1$.
By Lemma~\ref{le31} we have
\begin{eqnarray}\label{Z348}
M_{[\tau]}( Z) & = & M_{[S_1]^{\ep_1} \cdots [S_m]^{\ep_m}}( Z ) \nonumber \\
& = & M_{[S_1]^{\ep_1} \cdots [S_{m-1}]^{\ep_{m-1}}}(Z) + M_{[S_m]^{\ep_m}}( Z) 
\nonumber
~~~+~M_{[S_m]^{\ep_m}}( M_{[S_1]^{\ep_1} \cdots [S_{m-1}]^{\ep_{m-1}}}(  Z)) .
\end{eqnarray}
The induction hypothesis gives
\beql{Z349}
M_{[S_1]^{\ep_2} \cdots [S_{m-1}]^{\ep_{m-1}}}( Z) = \sum_{S \in \sG}
M_{[S]^{k'(S)}}(  Z)
\eeq
in which $k' (S) = k(S)$ if $S \neq S_m$ and $k'(S) = k(S) -\ep_m$ if $S=S_m$.
Applying this formula gives
\begin{eqnarray}\label{Z350}
M_{[S_m]^{\ep_m}} ( M_{[S_1]^{\ep_1} \cdots [S_{m-1}]^{\ep_{m-1}}}( Z)) 
& = &
\sum_{S \in \sG} M_{[S_m]^{\ep_m}}( M_{[S]^{k'(S)}}(Z)) \nonumber \\
& = & M_{[S_m]^{\ep_m}}(  M_{[S_m ]^{k'(S_m)}} ( Z)) \nonumber\\
& = & M_{[S_m]^{\ep_m + k' (S_m)}} (Z) - M_{[S_m]^{\ep_m}} (Z) -
M_{[S_m]^{k' (S_m)}}(Z), \nonumber
\end{eqnarray}
where we used \eqn{Z346} to remove all terms except $S=S_m$ and
applied Lemma \ref{le31} to obtain the last expression.
Thus
\begin{eqnarray*}
M_{[\tau ]} ( Z) & = & 
\left( \sum_{S \in \sG} M_{[S]^{k'(S)}}( Z) \right) + 
M_{[S_m]^{\ep_m}}( Z)  \\
&&
+ \,\,M_{[S_m]^{\ep_m + k' (S_m)}}(Z) - M_{[S_m]^{\ep_m}} (Z) - M_{[S_m]^{k' (S_m)}} (Z) \\
& = &
\left(\sum_{S \in \sG \atop S \neq S_m} M_{[S]^{k'(S)}}( Z )
 \right)
  + M _{[S_m]^{k' (S_m) +\ep_m}}( Z) 
  = \sum_{S \in \sG} M_{[S]^{k(S)}}( Z)~,
\end{eqnarray*}
which completes the induction step, and \eqn{Z343} follows.

The commutator subgroup of a free group on a set of generators $\sG$
is well known to be the set of words
$[\tau ] = [S_1]^{k_1} \cdots [S_m]^{k_m}$ in the generators
for which all exponents $k(S)$ satisfy
$$k(S):= \sum_{\{i:S_i =S\}} k_i =0 \quad\mbox{for all}\quad
S \in \sG ~.
$$
Now \eqn{Z343} gives
$$
M_{[\tau]} (Z) \equiv 0 \quad\mbox{for all}\quad [\tau] \in
\pi_1 (\sM , \bx_0)^{'} ~.
$$
The conclusion about $\tilde{\sM}^{ab}$ is immediate.
\end{proof}


\paragraph{\em Proof of Theorem \ref{th13}.}
 This is an
immediate consequence of Theorem~\ref{th31} and
Theorem~\ref{th32}. $\Box$\\

 \paragraph{\em Remark.} We have derived formulas
  for  $\sM$ as the base manifold,
rather than the extended manifold $\sM^{\#}$,  because
the largest manifold on which the analytically continued  functions $L^{\pm}(s, a, c)$
appearing in the functional
equation remain  well-defined (on a suitable covering manifold).
In addition  the manifold $\sM$ possesses a global automorphism of order $4$  corresponding
to the functional equation: $\phi_{\bR}: \sM \to \sM$, given by
$\phi_{\bR} (s, a, c) = (1-s, 1-c, a)$, which does not extend to $\sM^{\#}$.

%

\section{Differential-difference operators and Lerch monodromy functions}
\setcounter{equation}{0}

In the next three  sections we derive properties of the Lerch zeta
function and its monodromy functions $M_{[\tau]} (Z)$ as a function of the variables
$(s,a,c) \in \sM$. Here we show
that they  satisfy two differential-difference equations, and a  linear
partial differential equation in the operators
$\frac{\pt}{\pt a}$ and $\frac{\pt}{\pt c}$, on the universal cover $\tsM$.

The differential-difference equations involve the partial differential
operators
$\frac{\pt}{\pt a}$ and $\frac{\pt}{\pt c}$ viewed on $\sM$, and
extending to  the universal cover $\tilde{\sM}$. We define the {\em  raising operator}
\beql{N501}
\he{D}_L^{+} := \frac{\partial}{\partial c},
\eeq
and the {\em lowering operator}
\beql{N502}
\he{D}_L^{-} :=  \frac{1}{2 \pi i} \frac{\partial}{\partial a} + c.
\eeq
We define the {\em Lerch differential operator} as
\beql{Z46}
\he{D}_L := \he{D}_L^{-} \he{D}_L^{+} = \left( \frac{1}{2 \pi i} \frac{\pt}{\pt a} +c \right)
\frac{\pt}{\pt c} ~.
\eeq
The differential-difference equations satisfied by the Lerch zeta function in
the  region where the series expansion \eqn{101} converges  uniformly are
\beql{N503}
\he{D}_L^{+}(\zeta)(s, a, c) = -s \zeta(s+1, a ,c),
\eeq
and
\beql{N504}
\he{D}_L^{-}(\zeta)(s, a, c)= \zeta(s-1, a ,c).
\eeq
Combining \eqn{N502} and \eqn{N503} shows that the Lerch zeta function
satisfies the linear differential equation
$$
\he{D}_L \zeta(s, a, c) = -s  \zeta(s, a, c),
$$
a fact previously observed by Okubo \cite[p. 1057]{Ok98}.

Theorem~\ref{th42} below  gives the analytic continuation of these
equations. It uses the fact that the operators $\frac{\pt}{\pt a}$
and $\frac{\pt}{\pt c}$ lift to partial differential operators on
$\tsM$ which are {\em equivariant} with respect to the group
$\hG_{\sM}$ of diffeomorphisms of $\tilde{\sM}$ that commute with
the covering map $\pi: \tsM \to \sM$. Here
 $\hG_{\sM} \simeq \pi_1(\sM, \bx_0) \simeq \{ \sQ_{[\tau]} : [\tau ] \in \pi_1 (\sM ; \bf{x}_0 ) \}.$

%

\begin{theorem}\label{th42}
{\em (Lerch Differential-Difference Operators)}

(1) The analytic continuation $Z(s,a,c, [\gamma])$ of the Lerch zeta
function on the universal cover $\tsM$ satisfies the two differential-difference equations
\beql{Z46a}
\left( \frac{1}{2 \pi i} \frac{\pt}{\pt a} +c \right)Z(s,a,c, [\gamma]) =  Z (s-1,a,c, [\gamma_{-}])
\eeq
and
\beql{Z46b}
\frac{\pt}{\pt c} Z (s,a,c, [\gamma])= -s Z( s+1, a, c, [ \gamma_{+}]),
\eeq
in which $[\gamma_{+}]$ (resp. $[\gamma_{-}]$)  denote  paths which first traverse
$\gamma$ and then traverse a path from the endpoint of $\gamma$ that changes
the $s$-variable only, moving from  $s$ to $s + 1$ (resp. $s-1$).

(2) The analytic continuation $Z(s, a , c, [\gamma])$ on $\tsM$
satisfies the linear partial differential equation \beql{Z47}
\he{D}_L Z(s,a,c, [\gamma ]) = -sZ(s,a,c, [\gamma] ) ~, \eeq where
$\he{D}_L := \left( \frac{1}{2 \pi i} \frac{\pt}{\pt a}
\frac{\pt}{\pt c} +c \frac{\pt}{\pt c} \right).$

 (3) Each Lerch monodromy function $M_{[\tau]} (Z)(s, a, c, [\gamma])$
for each $[\tau ] \in \pi_1 ( \sM , \bx_0 )$ satisfies the two differential-difference
equations and the  differential equation above on  $\tsM$.
\end{theorem}

\begin{proof}
(1)  The Lerch zeta function $\zeta (s,a,c)$ is defined by
the Dirichlet series \eqn{101} in the simply-connected region
\beql{Z48}
 \sU = \{ s: s \in \CC\} \times \{ a: \Im (a) > 0 \} \times \{ c: \Re (c) > 0 \} ~. 
\eeq
In this region, $\zeta(s, a, c)$ satisfies the differential-difference equation
\begin{eqnarray}\label{Z49}
\left( \frac{1}{2 \pi i} \frac{\pt}{\pt a} +c \right)
\zeta (s,a,c) & = &
\frac{1}{2 \pi i} \sum_{n=0}^\infty (2 \pi in) e^{2 \pi ina} (n+c)^{-s} +
\sum_{n=0}^\infty ce^{2 \pi ina} (n+c)^{-s} \nonumber \\
& = & \zeta (s-1,a,c) ~.
\end{eqnarray}
It also satisfies the differential-difference equation
\begin{eqnarray}\label{Z410}
\frac{\pt}{\pt c} \zeta (s,a,c) & = & -s \sum_{n=0}^\infty e^{2 \pi ina} (n+c)^{-s-1} \nonumber \\
& = & -s \zeta (s+1,a,c) ~.
\end{eqnarray}
These equations carry over to $\tilde{\sM}$ by analytic continuation
using paths based at $\bx_0$. Note  that from any point $\bx= (s, a,
c, [\gamma]) \in \tsM$, there is a unique single-valued analytic
continuation possible in  the $s$-variable, in any  surface holding
the $(a, c)$-variables fixed. The difference operators $s \to s \pm
1$ are  interpreted to represent motions made in the $s$-variable on
such a  surface. The paths $[\gamma_{\pm}]$
represent the path $[\gamma]$ extended by a path in the $s$-variable
made in this surface, starting from the endpoint of $\gamma,$ moving
from $s$ to $s \pm 1$.

(2) Combining the differential-difference equations \eqn{Z49} and
\eqn{Z410} on the
 simply connected domain $\sU$   gives
 \beql{Z411}
\left( \frac{1}{2 \pi i} \frac{\pt}{\pt a} +c \right) \frac{\pt}{\pt c} \zeta (s,a,c) = -s \zeta (s,a,c) 
 \eeq 
valid on $\sU$. This shows
that \eqn{Z47} holds for $Z$ in the fundamental polycylinder
intersected with $\Im (a)> 0$. It now holds by analytic continuation
for all $Z([\gamma] )$ for all paths $\gamma$ based at $\bx_0$.

(3) By definition, we have
$$
M_{[\tau]} (Z) ( [\gamma]) :=
Z(s,a,c,[\tau \gamma ]) - Z(s,a,c,[\gamma] ) ~.
$$
Now we may apply the two differential-difference operators (and
their composition) $\he{D}_L^{\pm}$ to both sides of this equation.
By (1), (2) both terms on the right side satisfy these equations,
hence $M_{[\tau]} ([\gamma]) =M_{[\tau]} (Z) (s, a, c, [\gamma])$
satisfies \eqn{Z46a}, \eqn{Z46b}, \eqn{Z47} in place of 
$Z(s, a, c, [\gamma])$.
\end{proof}

\paragraph{{\bf Remarks.}}
(1) The relation $\he{D}_L f = -sf$ can be verified directly
for the monodromy formulae in Theorem \ref{th31},
and provides a consistency check on their correctness.

(2) The partial differential equation \eqn{Z47} is a  {\em local}
condition in a neighborhood of a point $(s,a,c ,[\gamma] ) \in \tsM$ and is
preserved under analytic continuation.
In contrast, the functional equation for the Lerch zeta function involves
difference operators,
which are not local.
Under analytic continuation equality is preserved in the functional
equation only if the proper
branches of the function are chosen, cf. Theorem~\ref{thK41} below.
For example the relation $\zeta (s,a,c) = \zeta (s,a+1,c)$ holds on the
domain \eqn{Z48}, but
under analytic continuation the right and left sides of this equality
may be on different sheets of the universal cover $\tsM.$

(3) All four terms in the functional equation in Theorem \ref{th11}
separately satisfy the partial differential equation of Theorem~\ref{th42}.
That is, \eqn{Z47} holds for $f(s,a,c)$ equal to any of
$\zeta (s,a,c)$, $e^{-2 \pi ia} \zeta (s,1-a,1-c)$, $e^{- 2 \pi iac} \zeta (1-s, 1-c, a)$ and
$e^{-2 \pi iac + 2 \pi ic} \zeta (1-s,c,1-a)$.
\vspace*{.2\baselineskip}

%

\section{Functional equations and monodromy functions}
\setcounter{equation}{0}

There are
linear relations between monodromy functions implied by
the functional equations. The  functional equation in Theorem~\ref{th11}
analytically
continues to the universal cover $\tsM$ of $\sM$,
provided each of the four terms continues along a separate path from the base
point $\bx_0 = (\frac{1}{2}, \frac{1}{2}, \frac{1}{2} )$.
Given a path $\gamma (t) = \{ (s(t), a(t), c(t)) : 0 \le t \le 1 \}$
in $\sM$ starting at
$\bx_0$,
we define the path
$\bth (\gamma )$ by
\beql{K41}
\bth (\gamma ) (t) :=
(1-s(t) , 1-c(t) , a(t)) ~,\quad 0 \le t \le 1 ~.
\eeq
In particular $\theta$ acts on loops in $\sM$ based at $\bx_0$.
This map on loops is well-defined on homotopy classes and induces a
homomorphism
\beql{K42}
\bth : \pi_1 (\sM , \bx_0) \to \pi_1 (\sM , \bx_0) ~,
\eeq
which is an automorphism of order 4.
%
\begin{theorem}\label{thK41}
(1) The map $\bth : \pi_1 ( \sM , \bx_0) \to \pi_1 ( \sM , \bx_0)$ is an
automorphism of order $4$ satisfying
\begin{eqnarray}\label{K43}
\bth ([X_n ]) & = & [Y_n ] \quad\mbox{for}\quad n \in \ZZ ~, \nonumber \\
\bth ([Y_n ]) & = & [X_{1-n} ] \quad\mbox{for}\quad n \in \ZZ ~.
\end{eqnarray}

(2) For each
$[\tau ] \in \pi_1 ( \sM, \bx_0)$
the Lerch zeta monodromy functions $M_{[\tau]} (Z)$ satisfy the
linear  relations
\begin{eqnarray}\label{K44}
&&\pi^{- \frac{s}{2}} \Gamma \left( \frac{s}{2} \right) \left\{
M_{[\tau]} (Z) ( [\gamma] ) + e^{- 2 \pi ia}
M_{\bth^2 ([\tau] )} (Z) ( [\bth^2 ( \gamma)] ) \right\} \nonumber \\
&=& e^{- 2 \pi iac} \pi^{- \frac{1-s}{2}} \Gamma \left(
\frac{1-s}{2} \right) \left\{ M_{\bth ([\tau])} (Z) ( [\bth ( \gamma
)]) + e^{2 \pi ic} M_{\bth^3 ([\tau])} (Z) ( [\bth^3 (\gamma )])
\right\}{~~~~~~~~~~}
\end{eqnarray}
and
\begin{eqnarray}\label{K45}
&& \pi^{- \frac{s+1}{2}}
\Gamma \left( \frac{s+1}{2} \right)
\left\{
M_{[\tau]} (Z) ( [\gamma] ) - e^{-2 \pi ia} M_{\bth^2 ([\tau ])}
(Z) ( [\bth^2 ( \gamma )]) \right\} \nonumber \\
&=& ie^{-2 \pi iac} \pi^{\frac{2-s}{2}} \Gamma \left( \frac{2 -
s}{2} \right) \left\{ M_{\bth ([\tau ])} (Z) ( [\bth ( \gamma )]) -
e^{2 \pi ic} M_{\bth^3 ([\tau ])} (Z) ([\bth^3 (\gamma )]) \right\}
{~~~~~~~~~~~}
\end{eqnarray}
in which $\gamma$ is any path with $\gamma (0)= \bx_0$.
\end{theorem}

\begin{proof}
The action \eqn{K43} of the automorphism $\bth$ on generators of
$\pi_1 (\sM , \bx_0)$ is easily computed from the definition. The
relations \eqn{K44} and \eqn{K45} now follow from analytically
continuing the functional equations \eqn{110b} and \eqn{111b} for
$\hat{L}^{\pm}(s,a,c)$
along a path \eqn{K41}.
\end{proof}

\paragraph{\em Remark.}
Theorem~\ref{thK41} provides another consistency check on the monodromy
functions. As an example, taking
$[\tau] = [X_0]$ and using Theorem~\ref{th31} to compute monodromy functions
one finds that \eqn{K44} asserts that
\begin{eqnarray}\label{K46}
&&\pi^{- \frac{s}{2}}\Gamma \left( \frac{s}{2} \right)
\left\{ - \frac {(2 \pi)^s e^{ \frac{\pi i s}{2}}}{\Gamma (s)}
a^{s - 1} e^{ - 2 \pi i a c} -
 e^{- 2 \pi i a}  \frac {(2 \pi)^s e^{ \frac{\pi i s}{2}}}{\Gamma (s)}
e^{\pi i (s - 1)}a^{s - 1} e^{-2 \pi i (1-a)(1-c) + 2 \pi i (1 - c)} \right\}    \nonumber \\
&&~~~~ = e^{- 2 \pi i a c} \pi^{- \frac{1-s}{2}}
\Gamma \left( \frac{1-s}{2} \right)
\left\{ (e^{2 \pi i s} - 1) a^{-(1 - s)}  +~~~ 0 \right\}.
\end{eqnarray}
This can be verified directly using suitable gamma function identities.

%

\section{Lerch monodromy vector spaces}
\setcounter{equation}{0}

In this section we determine properties of
the restricted monodromy functions obtained when the parameter
$s$ is viewed as fixed. The values $s \in \ZZ$ are
``special values,'' where the vector space spanned by the
monodromy functions simplifies.

For fixed $s \in \CC$ and $[\tau] \in \pi_1(\sM, \bx_0)$, denote by
$\sQ_{\tau}^s(Z)$ (resp.  $M_{[\tau]}^s (Z)$)  the restriction of the functions $\sQ_{\tau}(Z)$ (resp.
$M_{[\tau]} (Z)$) to the manifold 
 \beql{Z41} 
 \tsM_s := \pi^{-1}(\{
s \} \times (\CC \smallsetminus \ZZ) \times( \CC \smallsetminus
\ZZ)) 
\eeq 
which is a path-connected analytic submanifold of the universal cover $\tsM$ of $\sM$.
%
%
\begin{defi}
{\em The  {\em Lerch monodromy space}
$\sV_s$ {\em at} $s \in \CC$ is  the 
vector space generated by all  germs of function elements of $Z$, which is generated by
$Z^s$ and by all
$\{ Q_{[\tau]}^s(Z):~ [\tau] \in \pi_1(\sM, \bx_0)\},$
where $Q_{[\tau]}^s(Z)$
 is the function $\sQ_{[\tau]}(Z)$ restricted to $s \in \CC$ as above.
This vector
space consists of all finite $\CC$-linear combinations of the specified
generators $\sQ_{[\tau]}^s(Z)$.}
\end{defi}

It is immediate that this vector space is spanned by 
the function $Z^s =\sQ_{[id]}^s(Z)$ together with all the functions $M_{[\tau]}^s (Z)$ regarded as
functions on $\tsM_s$,
 i.e.
 \beql{Z42}
 \sV_s := \CC[Z] + \left( \sum_{[\tau]\in \pi_1 (\sM , \bx_0)} \CC[M_{[\tau]}^s (Z)] \right)~,
\eeq 
using the fact that $M_{[\tau]}^s(Z) = \sQ_{[\tau]}^s (Z) - Z^s.$

%
%

\begin{theorem}\label{th41}
{\em (Lerch Monodromy Spaces)}
The Lerch monodromy space $\sV_s$ depends on the parameter $s \in \CC$
as follows.
\begin{itemize}
\item[(i)]
If $s \in \ZZ$ with $s = - m \le 0$, then all Lerch monodromy functions
vanish identically, i.e.
\beql{Z43}
M_{[\tau]}^{-m} (Z) =0 \quad\mbox{for all}\quad [\tau ] \in \pi_1
(\sM , \bx_0 ) ~.
\eeq
Thus $\sV_{-m} = \CC[Z^{-m}]$ is one-dimensional.
\item[(ii)]
If $s \in \ZZ$ with $s= m \ge 1$, then $\sV_m$ is an infinite-dimensional
vector space,
and has a basis consisting of the function $Z^m$ together with 
$\{M_{[X_n]}^m (Z) : n \in \ZZ \}$.
\item[(iii)]
If $s \not\in \ZZ$, then $\sV_s$ is an infinite-dimensional vector space,
and has as a basis the function $Z^s$ together with the
functions
$\{M_{[X_n]}^s (Z) : n \in \ZZ \} \cup \{ M_{[Y_n]}^s :
n \in \ZZ_{\le 0} \}$.
\end{itemize}
\end{theorem}

\begin{proof}
(i)~Let $s=-m \in \ZZ_{\le 0}$. Theorem \ref{th31} shows that the
monodromy $M_{[\tau]}^s (Z)$ vanishes identically when 
$[\tau ] =[S]^k$ for any generator $S \in \sG$. Indeed, this holds for 
$[\tau]= [X_n]^k$ because $\frac{1}{\Gamma (s)} =0$ when 
$s \in \ZZ_{\le 0}$ and it holds for $[\tau] = [Y_n]^k$ 
because $e^{-2 \pi is} -1=0$ for $s \in \ZZ$. Now Theorem \ref{th32} shows that all monodromy
functions $M_{[\tau]}^{-m}(Z)$ vanish identically. 
In this case \eqn{Z42} shows that $\sV_{-m} = \CC[Z^m]$ is one-dimensional.

(ii)~For $s=m \in \ZZ_{>0}$ the monodromy functions $M_{[Y_n]^k}^s (Z)$
vanish identically as in part  (i). However
Theorem \ref{th31} shows that for $n \in \ZZ$ the monodromy function
\[
M_{[X_n]}^s (Z)= c_n (a-n)^{s-1} e^{-2 \pi i(a-n)c},
\]
in which
$c_n$ is a {\em nonzero} constant (which itself depends on $s$).
Furthermore all $M_{[X_n]^k}^s (Z)$ have exactly the same form, so are
linearly dependent on $M_{[X_n]}^s (Z)$.
Theorem \ref{th32} now shows that $\{ M_{[X_n]}^m (Z): n \in \ZZ \}$ together with $Z^m$
form a spanning set for the vector space
$\sV_m$.
It is a basis because any finite subset is linearly independent,
by considering them as functions of $(a,c)$ in a small disk around
$(s, \frac{1}{2} , \frac{1}{2})$ with $s$ held constant.
Thus $\sV_m$ is infinite-dimensional.

(iii)~By a similar argument to (ii), the functions in
$$
\sB := \{ M_{[X_n]}^s (Z) : n \in \ZZ \}
\cup \{ M_{[Y_n]}^s (Z) : n \in \ZZ_{\le 0} \} \cup \{Z^s\}
$$
comprise a spanning set of the vector space $\sV_s$.
The functions $M_{[Y_n]}^s (Z)$ have the form
$c_n e^{-2 \pi ina} (c-n)^{-s}$ for {\em nonzero}
$c_n$ (which depends on $s$).
The $\CC$-linear independence of any finite subset of the functions
in $\sB$ can be established by considering them in a small disk
in the $(a,c)$-plane around $(s, \frac{1}{2} , \frac{1}{2} )$, holding $s$
constant. Thus $\sV_s$ is infinite-dimensional.
\end{proof}

We can associate to the vector space $\sV_s$ an induced {\em monodromy representation}
\beql{Z44}
\rho_s : \pi_1 ( \sM , \bx_0) \to End (\sV_s ) ~,
\eeq
defined for $[\sigma] \in  \pi_1 ( \sM , \bx_0)$ by
\beql{Z45}
\rho_s ([\sigma]) (\sQ_{[\tau]}^s (Z)) := \sQ_{[\sigma^{-1} \tau]}^s (Z).
\eeq
Now one has
\beql{745a}
\rho_s([\sigma_1])  \rho_s([\sigma_2]) = \rho_s([ \sigma_1 \sigma_2]),
\eeq 
and we note that 
\begin{eqnarray}\label{Z45b}
\rho_s([\sigma] ) (M_{[\tau]}^s (Z)) &= & \rho_s([\sigma]) \left(\sQ_{[\tau]}^s(Z) - \sQ_{\mbox{Id}}^s(Z)\right)\nonumber \\
&=& \sQ_{[\sigma^{-1} \tau]}^s(Z) - \sQ_{[\sigma^{-1}]}^s(Z) \nonumber \\
&=& M_{[\sigma^{-1} \tau]}^s(Z) - M_{[\sigma^{-1}]}^s(Z).
\end{eqnarray}
This representation has a kernel that depends on $s$,
since there are relations among $M_{[S]^d}^s (Z)$
for a fixed generator $[S]$ as $d \in \ZZ$ varies.
For $s \not\in \ZZ$ the quotient group
$$
\hG_s := \pi (\sM, \bx_0) / \ker (\rho_s )
$$
 can be identified with the free
abelian group $\ZZ [ \tilde{\sG} ]$ on the
generating set
$\tilde{\sG}= \{ [X_n] : n \in \ZZ \} \cup \{ [Y_n] : n \le 0 \}$,
and the quotient representation
$\bar{\rho}_s : \ZZ [\tilde{\sG} ] \to End ( \sV_s )$ studied as $s$ varies.\\

\paragraph{\em Remarks.}
(1) The multi-valued nature of the Lerch zeta function
encodes information about the interaction of the additive and
multiplicative structures of $\ZZ$.
The vanishing of all monodromy functions
when $s=-m$ is a  non-positive integer  is another aspect of
the viewpoint that these
points are  ``special values'' that contain ``universal'' information.

(2) In part III \cite[Sect. 5]{LL3} we will show  that for such
$s=-m$, taking ``special values'' $L^{\pm}(-m, a, c)$ for $n \in
\ZZ_{>0}$ at  rational values of $a$ and $c$ with $0 < a, c \leq 1$
yields data sufficient to construct by $p$-adic interpolation all
$p$-adic $L$-functions associated to $\QQ$.

%

\section{Extended analytic continuation}
\setcounter{equation}{0}

In this section we complete the analytic continuation of the Lerch
zeta function $\zeta(s, a, c)$ to a multivalued function $Z(s, a, c,
[\gamma])$ defined over the manifold \beql{801} \sM^{\#} := \{ (s,
a, c) \in \CC \times (\CC \smallsetminus \ZZ) \times (\CC
\smallsetminus \ZZ_{\le 0})\}, \eeq obtained from $\sM$ by filling
in all points $\{ c: c=n \ge 1\}$. The  vanishing of the Lerch
monodromy functions \beql{802} M_{[Y_n]} (Z) =0 \quad\mbox{for}\quad
n=1,2,3, \ldots  \eeq (as given in  Theorem~\ref{th31}) implies that
the function $Z(s, a, c, [\gamma])$ is single-valued in a
(``punctured") open neighborhood of  any point $(s, a, c) \in \sM$,
with $c = n$ for integer $n \ge 1$. However a proof is required that
the singularity at $c=n$ can be removed.

\paragraph{\em Proof of Theorem~\ref{th23}. }
The proof parallels the analytic continuation to $\sM$, so we omit
many details. We must show that for each $(s, a) \in \CC \times (\CC
\smallsetminus \ZZ)$ and for each $[\gamma]$ with endpoint $\bx_1$
lying over $(s, a, n) \in \sM^{\#}$, there is a local  open
neighborhood  in three complex variables where the function $Z(s,a ,
n, [\gamma])$ is holomorphic. This reduces to showing holomorphicity
in three variables for $\zeta(s, a, c)$ on the principal sheet, plus
holomorphicity  of each possible monodromy function $M_{[\tau]}(s,
a, c)$ on the principal sheet. The form of the monodromy functions
in Theorem~\ref{th31} and Theorem~\ref{th32} indicates that all of
them are holomorphic in an open three-dimensional neighborhood of
any point with $c=n$, so it remains to check the case of $\zeta(s,
a, c)$ itself.

 First, the earlier proof showed that  Dirichlet series expansion \eqn{101} gives
an analytic continuation of $\zeta(s, a, c)$ to the region
$$
 \sU  = \{ s:  s \in \CC \} \times \{ a: \Im (a) > 0 \}
\times \{c: \Re (c) > 0 \} ~. 
$$
Second, the integral
representation 
\beql{804}
 \zeta (s,a,c) := \frac{1}{\Gamma (s)} 
 \int_0^\infty \frac{t^{s-1} e^{-ct}}{1- e^{2 \pi ia} e^{-t}} dt ~.
\eeq 
defines  a single-valued  analytic extension to the region
\beql{805}
 \sW_1= \{s: \Re (s) > 0 \} \times \{a: 0< \Re(a) < 1 \}
\times \{c: \Re (c) > 0 \} ~. 
\eeq 
Third, the periodicity of the
function $\zeta(s, a, c)$ under $a \mapsto a+1$ in the region 
$\sU$,
when compared to $\sW_1$,  means that we can use this periodicity to
extend the continuation in $\sW_1$ to the lower half-plane in strips
of width one, to the region
 \beql{806} \sW_2 = \{s: \Re (s) > 0 \}
\times \{a: a \in \sA_L \} \times \{c: \Re (c) > 0 \} 
~ \eeq
 in which
$$
\sA_L := \CC \smallsetminus \{ M_k : k \in \ZZ \} ~
$$
is pictured in Figure \ref{fg31}. Fourth, using the monodromy
functions $M_{[X_n]}(Z)$ to circle around all integer points in the
$a$-plane, we can now analytically continue  this function to a
multi-valued function over the manifold
$$
\sW_3= \{s: \Re (s) > 0 \} \times \{a: a \in \CC \smallsetminus \ZZ \} \times \{c: \Re (c) > 0 \},
$$
which will be  single-valued on the maximal abelian cover of this manfold.

Fifth, we use the differential-difference equation \beql{807}
\zeta(s-1, a, c) =\left(\frac{1}{2\pi i} \frac{\partial}{\partial a}
+ c\right)(\zeta)(s, a, c) \eeq to analytically continue from $s$ to
$s-1$ in the $s$-variable, from the region $\sW_3$. Here we use
Theorem~\ref{th42} and the fact that the differential-difference
equation continues to hold at points where  $c=n \ge 1$, in the
region $\sW_3$, by analytic continuation. We now have obtained
analytic continuation as a multi-valued function over the manifold
$$
\sW_4 = \{s: s \in \CC \} \times \{a: a \in \CC \smallsetminus \ZZ \} \times \{c: \Re (c) > 0 \},
$$
which is single-valued over the maximal abelian cover.
(Note that we cannot use the functional equation  to make the continuation
in $s$, at points with $c=n \ge 1$.)

Sixth, we extend the analytic continuation to the rest of the  $c$-plane,
omitting non-positive integers, using
Theorem~\ref{th13} and the
fact that this function agrees with $Z(s, a, c, [\gamma])$ in the resulting region.
Since all resulting  monodromy functions are holomorphic at points where $c= n \ge 1$,
we obtain a multi-valued function defined over the manifold
$$
\sM^{\#} =\{ (s, a,c) \in  \CC \times \left( \CC \smallsetminus
\ZZ\right) \times \left(\CC \smallsetminus \ZZ_{\le 0}\right)\}.
$$
Since no  monodromy is added at points $c=n$, the fact
that this function is single-valued on the maximal
abelian cover $\tsM^{\#}$ of $\sM^{\#}$
 follows from Theorem~\ref{th13}.
$~~~\bsq$\\

\paragraph{\em Remarks.}
(1) Applying the  transformation $\phi_{\bR}(s, a,c) = (1-s, 1-c,
a)$ defined on $\CC^3$,
 we deduce that the function
$\he{R}(\zeta)(s, a, c):= e^{2\pi i ac}\zeta(s, 1-c, a)$ implicitly appearing
in the extended Lerch transformation formula  (\ref{224a}) has a multi-valued analytic
continuation defined over the extended manifold
 \beql{809}
\phi_{\bR}(\sM^{\#}) = \CC \times \left( \CC \smallsetminus \ZZ_{\ge 1} \right)
 \times \left(\CC \smallsetminus \ZZ \right).
 \eeq
Similar results hold for $\he{R}^k(\zeta)$ for $k = 2, 3$, giving
analytic continuations above the manifolds $\phi_{\bR}^k(\sM^{\#})$.

(2) The proof above used the differential-difference operator \eqn{N502} to
analytically continue from $\Re(s) >0$ to $s \in \CC$. An
alternative approach to analytically continue in the $s$-variable is
to use a different contour  integral representation of the Lerch
zeta function, as done by  Barnes \cite{Ba06}. That integral  uses
the same right hand side in \eqn{804}, with the integral  taken over
a keyhole contour on the complex plane cut along the positive real
axis, that initially comes in from $+\infty$ along the (lower)
positive real axis, circles clockwise around $s=0$ on a small circle
and goes back to $+ \infty$ along the (upper) real axis. This
integral equals $(e^{2 \pi is} -1) \zeta(s,a, c)$.

%
\section{Concluding remarks}
\setcounter{equation}{0}

The  analytic continuation of the Lerch zeta function given in this
paper is to an essentially maximal domain of holomorphy in
three complex variables.
It  is also possible to  define
analytic continuations of the Lerch zeta function in fewer variables which make sense on
various  ``singular strata", i.e. regions $(s, a,c)$ in which 
either $a$ or $c$, or both, take integer values, corresponding to the punctures in
the manifold $\sM^{\#}$.  These functions can sometimes be defined as continuous
limits from ``non-singular"  ranges of $(s, a,c)$, as discussed in part I, and are then analytic in
fewer variables on such ranges.  
For example, the  Hurwitz zeta function corresponds to fixing $a=0$ and taking $0<c<1$,  and 
we view it  as belonging to 
a  ``singular stratum" of  real codimension two (i.e. this stratum is a $2$-dimensional
complex manifold, where we vary the other two variables.)
As noted in part I,  it can be obtained 
as a limiting value of the Lerch zeta function on $\sM$, provided that $\Re(s) > 1$. 
It is well-known that the Hurwitz zeta function can be 
meromorphically continued in the $s$-variable, and has a simple
pole at $s=1$. Our viewpoint here is that the polar singularity is a signal 
that this function  belongs to a ``singular stratum."  One can also show by the methods of this
paper that it  
possesses a multi-valued
analytic continuation in the $(s,c)$-variables, having 
additional singularities at non-positive integer
values of $c$.
The  Riemann zeta function is associated to various integer
points in both the $a$-variable and $c$-variable. For the point  $a=1, c=1$ it
lies on a singular stratum of real codimension $2$, and arises as a continuous limiting
value of the Lerch zeta function on $\sM$, when $\Re(s) > 1$. 
According to part I, it  is also (in some sense) associated to  the point 
$a=0, c=0$, which  lies on a singular stratum of real codimension $4$,
but there the Lerch zeta  function has no continuous limiting value
for any $s \in \CC$. We note that 
specialization of the  Lerch functional equation in Theorem~\ref{th11}
at $(a,c)=(1,1)$ necessarily requires that data at  $(a, c)= (0,0)$ be included.
 It seems  an interesting topic to further identify and classify all
the ``singular strata" and to understand 
 how the functions (or collection of functions)  on each singular stratum
arise  by ``degeneration" from the multivalued function above.

Finally we have observed  that  the strata corresponding to positive integer $c= n \ge 1$ 
are in fact  ``nonsingular," as indicated by Theorem \ref{th23}. The truly  ``singular strata"  can be
recognized by the appearance of discontinuities, already illustrated
in part I \cite[Theorem 2.3]{LL1}.  


%


\begin{thebibliography}{99}
\bibitem{Ap51}
T. M. Apostol,
On the Lerch zeta function,
{\em Pacific J. Math.} {\bf 1} (1951), 161--167.

\bibitem{Ap76}
T. M. Apostol,
{\em Introduction to Analytic Number Theory}, Springer-Verlag:
New York, 1976.

\bibitem{Ba06}
E. W. Barnes,
On certain functions defined by Taylor's series with
finite radius of convergence,
Proc. Lond. Math. Soc. (Ser. 2) {\bf 4} (1906), 287--316.


\bibitem{Be72}
B. Berndt,
Two new proofs of Lerch's functional equation,
{\em Proc. Amer. Math. Soc.} {\bf 32} (1972), 403--408.

\bibitem{Er53}
A. Erd\'{e}lyi et~al.,
{\em Higher Transcendental Functions -- Vol. I},
M. Graw-Hill: New York 1953.

\bibitem{Ha01}
A. Hatcher,
{\em Algebraic Topology}, 
Cambridge University Press,
Cambridge 2001.

\bibitem{Hi76}
E. Hille,
{\em Ordinary Differential Equations in the Complex Domain,}
John Wiley \& Sons:  New York 1976.
[Reprint: Dover Publications: New York 1997.]


\bibitem{Hu82}
A. Hurwitz, Einige Eigenshaften der Dirichlet'schen
Funktionen $F(s) = \sum \left( \frac{D}{n} \right) \frac{1}{n^s}$,
dei bei der Bestimmung der Klessenzahlen binarer
quadratischen formen auftreten,
{\em Zeitschrift f\"{u}r Math. und Physik} {\bf 27} (1882), 86--101.

\bibitem{KKY00}
S. Kanemitsu, M. Katsurada and M. Yoshimoto,
On the Hurwitz-Lerch zeta-function,
Aequationes Math. {\bf 59}(2000), 1--19.

\bibitem{KT07}
S. Kanemitsu and H. Terada,
{\em Vistas of Special Functions},
World Scientific: Singapore 2007.





\bibitem{LL1}
J. C. Lagarias and W-C. W. Li,
The Lerch zeta function I. Zeta integrals,
preprint.

\bibitem{LL3}
J. C. Lagarias and W-C. W. Li,
The Lerch zeta function III. Polylogarithms and special values,
preprint.

\bibitem{LL4}
J. C. Lagarias and W.-C. W. Li,
The Lerch Zeta Function IV. Hecke Operators,
in preparation.

\bibitem{LG02}
A. Lauren\v{c}ikas and R. Garunk\v{s}tis,
{\em The Lerch zeta-function},
Kluwer Academic Publishers: Dordrecht 2002.

\bibitem{Le87}
M. Lerch,
Note sur la fonction
$\sR (w,x,s) = \sum_{k=0}^\infty \frac{e^{2k \pi ix}}{(w+k)^s}$,
{\em Acta Math.} {\bf 11} (1887), 19--24.


\bibitem{L2}
L. Lewin, Ed.,
{\em Structural Properties of Polylogarithms},
AMS: Providence, RI 1991.

\bibitem{Li57}
M. Lipschitz,
Untersuchung einer aus vier Elementen gebildeten Reihe,
{\em J. Reine Angew.} {\bf 54} (1857), 313--328.

\bibitem{Li89}
M. Lipschitz,
Untersuchung der Eigenshaften einer Gattung von unendlichen Reihen,
{\em J. Reine Angew.} {\bf 105} (1889), 127--156.

\bibitem{Mi71}
M. Mikolas,
New proof and extension of the functional equation of Lerch's
zeta function,
{\em Ann. Univ. Sci. Budapest - E\"{o}tvos Sect. Math.} {\bf 14} (1971),
111--116.


\bibitem{Mu83}
D. Mumford,
{\em Tata Lectures on Theta I},
Birkh\"{a}user: Boston 1983.

\bibitem{Ob56}
E. Oberhettinger,
Note on the Lerch zeta function, {\em Pacific J. Math.}
{\bf 6} (1956), 117--120.

\bibitem{Ok98}
S. Okubo,
Lorenz-Invariant Hamiltonian and Riemann Hypothesis,
{\em J. Phys. A} {\bf 31} (1998), 1049--1057.

\bibitem{Pa88}
S. J. Patterson,
{\em An introduction to the theory of the Riemann zeta function},
Cambridge University Press, Cambridge 1988.


\bibitem{SC01}
H. M. Srivastava and J. Choi,
{\em Series Associated with the Zeta and Related Functions},
Kluwer Academic Publishers: Dordrecht 2001. 


\bibitem{Ta67}
J. Tate,
Fourier Analysis in Number Fields and Hecke's Zeta-Functions,
in: {\em Algebraic Number Theory} (J. W. S. Cassels and A. Fr\"{o}hlich, Eds.),
Academic Press: London 1967, pp. 305--347.

\bibitem{We76}
A. Weil,
{\em Elliptic Functions according to Eisenstein and Kronecker},
Springer-Verlag, New York 1976.

\bibitem{WW27}
E. T. Whittaker and G. N. Watson,
{\em A Course of Modern Analysis},
Fourth Edition, 1927.
Cambridge University Press: Cambridge: Reprint 1965.


\end{thebibliography}
\end{document}